\documentclass{amsart}

\newtheorem{theorem}{Theorem}[section]

\newtheorem{cor}[theorem]{Corollary}

\theoremstyle{definition}
\newtheorem{definition}[theorem]{Definition}

\theoremstyle{remark}

\numberwithin{equation}{section}

\usepackage{amsmath}
\usepackage{amsfonts}
\usepackage{amssymb}

\newcommand{\abs}[1]{\ensuremath{\left| #1 \right| }}
\newcommand{\qPs}{$q$--Pochhammer symbol}

\newcommand{\rhs}{right hand side}
\newcommand{\lhs}{left hand side}

\newcommand{\wrt}{with respect to}

\newcommand{\Wt}{Watson transformation}

\newcommand{\RRis}{Rogers--Ramanujan identities}

\newcommand{\Js}{Jackson sum}

\newcommand{\ci}{cocycle identity}

\newcommand{\tns}{\ensuremath{_{10}\varphi_9}}

\newcommand{\Jc}{Jackson coefficients}
\newcommand{\eJc}{elliptic Jackson coefficients}
\newcommand{\Mf}{Macdonald functions}

\begin{document}

\title[Multiple basic special numbers]{Multiple analogues of binomial coefficients and related families of special numbers} 

\author{Hasan Coskun}
\address{Department of Mathematics, Texas A\&M
  University--Commerce, Binnion Hall, Room 314, Commerce, TX 75429}  
\curraddr{Department of Mathematics, Texas A\&M
  University--Commerce, Binnion Hall, Room 314, Commerce, TX 75429}  
\email{hasan\_coskun@tamu-commerce.edu}

\date{December 21, 2009} % and, in revised form, June 22, 1994.}


\begin{abstract}
We construct multiple $qt$-binomial coefficients and related multiple analogues of several celebrated families of special numbers in this paper. These  multidimensional generalizations include the first and the second kind of $qt$-Stirling numbers, $qt$-Bell numbers, $qt$-Bernoulli numbers, $qt$-Catalan numbers and the $qt$--Fibonacci numbers.   
In the course of developing main properties of these extensions, we prove results that are significant in their own rights such as certain probability measures on the set of integer partitions.
\end{abstract}

\maketitle

\section{Introduction}
\label{introduction}
Many distinct sequences of special numbers are investigated and their properties are explored in number theory.
Various generalizations of such families of numbers have also been studied extensively for most cases. An important class of generalizations of special numbers is their one dimensional, one parameter $q$-extensions. In this paper, we give multidimensional basic $qt$-generalizations for several collections of such numbers including the binomial coefficients and the classes of numbers associated with the names of Stirling, Bernoulli, Catalan, Bell and Fibonacci. We also point out how to construct multiple ordinary $\alpha$-generalizations for the same number sequences.  
The definitions of these numbers and the properties they satisfy show great variety which makes this research area very interesting.  

Among many mathematicians who contributed to this line of research, L. Carlitz appears to be the first to study the $q$-extensions for several families of special numbers given in this paper. Many have made significant contributions since then investigating properties of $q$-generalizations, their applications and connections with other types of numbers. We will give references to some of these successful efforts and important work in section~\ref{mulipleSpecial} below. 

The present paper takes a step in generalizing the one dimensional $q$-special numbers to multiple $qt$-special analogues. These generalizations are given in terms of the $qt$--binomial coefficient defined in section~\ref{section1} in the form
\begin{equation*}
\binom{\lambda}{\mu}_{\!\!\!q,t} := \dfrac{ q^{|\mu|} t^{2n(\mu)+(1-n)|\mu| } }  { (qt^{n-1} )_\mu} \prod_{1\leq i<j \leq n} \left\{\dfrac{ (qt^{j-i})_{\mu_i-\mu_j} } {(qt^{j-i-1})_{\mu_i-\mu_j}  } \right\}  W^{s\uparrow}_\mu(q^\lambda t^{\delta(n)}; q, t)
\end{equation*}
where $\lambda$ and $\mu$ are $n$-part partitions and $q,t\in\mathbb{C}$. The $W^{s\uparrow}_\mu$ function that enters into the definition is a limiting case of the well--poised $BC_n$ Macdonald function $W_\lambda$. We first give a brief review of this remarkable family of functions and its one parameter generalization $\omega_{\lambda}$ called well--poised $BC_n$ Jackson coefficients.   % that we need in the proof of our results. 
%This will be followed by a discussion the basic one parameter $BC_n$ \BL\ proved in~\cite{Coskun2}. 
The symmetric rational functions $W_\lambda$ and $\omega_{\lambda}$ are first introduced in the author's Ph.D. thesis~\cite{Coskun0} in the basic (trigonometric) case, and later in~\cite{CoskunG1} in the more general elliptic form.   

\section{Background}

We start with the definition of the \qPs\ $(a;q)_\alpha$,
where $q, \alpha\in\mathbb{C}$, which can be defined formally by
\begin{equation}
\label{qPochSymbol}
(a;q)_\alpha :=\dfrac{(a;q)_\infty}{(aq^\alpha;q)_\infty}
\end{equation} 
in terms of the infinite product $(a;q)_\infty:=\prod_{i=0}^{\infty} (1-aq^i)$. Note that when $\alpha=m$ is a positive integer, the definition reduces to $(a;q)_m= \prod_{k=0}^{m-1}(1-aq^k)$.

An elliptic
analogue is defined~\cite{FrenkelT1, RosengrenS1} in the form 
\begin{equation}
  (a; q,p)_m := \prod_{k=0}^{m-1} \theta(aq^m)
\end{equation}
where $a\in \mathbb{C}$, $m$ is a positive integer and the normalized elliptic
function $\theta(x)$ is given by
\begin{equation}
  \theta(x) = \theta(x;p) := (x; p)_\infty (p/x; p)_\infty
\end{equation}
for $x, p\in \mathbb{C}$ with $\abs{p}<1$. The definition is extended to negative $m$ by setting $(a; q,p)_m = 1/ (aq^{m}; q, p)_{-m} $.  Note that when $p=0$, $(a; q,p)_m$ reduces to the standard (trigonometric) \qPs.

For any partition $\lambda = (\lambda_1, \ldots, \lambda_n)$ and
$t\in\mathbb{C}$, define~\cite{Warnaar2}
\begin{equation}
\label{ellipticQtPocSymbol}
  (a)_\lambda=(a; q, p, t)_\lambda := \prod_{k=1}^{n} (at^{1-i};
  q,p)_{\lambda_i} .
\end{equation}
Note that when $\lambda=(\lambda_1) = \lambda_1$ is a single part
partition, then $(a; q, p, t)_\lambda = (a; q, p)_{\lambda_1} =
(a)_{\lambda_1}$. The following notation will also be used.
\begin{equation}
  (a_1, \ldots, a_k)_\lambda = (a_1, \ldots, a_k; q, p, t)_\lambda :=
  (a_1)_\lambda \ldots (a_k)_\lambda .
\end{equation}

Now let $\lambda=(\lambda_1, \ldots, 
\lambda_n)$ and $\mu=(\mu_1, \ldots, \mu_n)$ be partitions of at most
$n$ parts for a positive integer $n$ such that the
skew partition $\lambda/\mu$ is a horizontal strip; i.e. $\lambda_1
\geq \mu_1 \geq\lambda_2 \geq \mu_2 \geq \ldots \lambda_n \geq
\mu_n \geq \lambda_{n+1} = \mu_{n+1} = 0$. Following~\cite{Coskun1}, we define
\begin{multline}
\label{definitionHfactor}
H_{\lambda/\mu}(q,p,t,b) \\
:= \prod_{1\leq i < j\leq
n}\left\{\dfrac{(q^{\mu_i-\mu_{j-1}}t^{j-i})_{\mu_{j-1}-\lambda_j}
(q^{\lambda_i+\lambda_j}t^{3-j-i}b)_{\mu_{j-1}-\lambda_j}}
{(q^{\mu_i-\mu_{j-1}+1}t^{j-i-1})_{\mu_{j-1}-\lambda_j}(q^{\lambda_i
    +\lambda_j+1}t^{2-j-i}b)_{\mu_{j-1}-\lambda_j}}\right.\\
\left.\cdot 
\dfrac{(q^{\lambda_i-\mu_{j-1}+1}t^{j-i-1})_{\mu_{j-1}-\lambda_j}}
{(q^{\lambda_i-\mu_{j-1}}t^{j-i})_{\mu_{j-1}-\lambda_j}}\right\}\cdot\prod_{1\leq
i <(j-1)\leq n}
\dfrac{(q^{\mu_i+\lambda_j+1}t^{1-j-i}b)_{\mu_{j-1}-\lambda_j}}
{(q^{\mu_i+\lambda_j}t^{2-j-i}b)_{\mu_{j-1}-\lambda_j}}
\end{multline}
and 
\begin{multline}
\label{definitionSkewW}
W_{\lambda/\mu}(x; q,p,t,a,b)
:= H_{\lambda/\mu}(q,p,t,b)\cdot\dfrac{(x^{-1}, ax)_\lambda
  (qbx/t, qb/(axt))_\mu}
{(x^{-1}, ax)_\mu (qbx, qb/(ax))_\lambda}\\
\cdot\prod_{i=1}^n\left\{\dfrac{\theta(bt^{1-2i}q^{2\mu_i})}{\theta(bt^{1-2i})}
  \dfrac{(bt^{1-2i})_{\mu_i+\lambda_{i+1}}}
{(bqt^{-2i})_{\mu_i+\lambda_{i+1}}}\cdot
t^{i(\mu_i-\lambda_{i+1})}\right\}
\end{multline}
where $q,p,t,x,a,b\in\mathbb{C}$. The function
$W_{\lambda/\mu}(y, z_1, \ldots, 
z_\ell; q,p,t,a,b)$ is extended to $\ell+1$ variables $y, z_1, \ldots, z_\ell
\in\mathbb{C}$  
through the following recursion formula
\begin{multline}
\label{eqWrecurrence}
W_{\lambda/\mu}(y,z_1,z_2,\ldots,z_\ell;q, p, t, a, b) \\
= \sum_{\nu\prec \lambda} W_{\lambda/\nu}(yt^{-\ell};q, p, t, at^{2\ell},
bt^\ell) \, W_{\nu/\mu}(z_1,\ldots, z_\ell;q, p, t, a, b).
\end{multline}

We will also need the \eJc\ below. Let $\lambda$ and
$\mu$ be again 
partitions of at most $n$--parts such that $\lambda/\mu$ is a skew
partition. Then the \Jc\ $\omega_{\lambda/\mu}$ are defined by
\begin{multline}
\label{eq:omega{lambda,mu}}
\omega_{\lambda/\mu}(x; r, q,p,t; a,b)
:= \dfrac{(x^{-1}, ax)_{\lambda}} {(qbx, qb/ax)_{\lambda}}
    \dfrac{(qbr^{-1}x, qb/axr)_{\mu}}{(x^{-1}, ax)_{\mu}} \\
\cdot \dfrac{(r, br^{-1}t^{1-n})_{\mu}}{(qbr^{-2}, qt^{n-1})_{\mu}}
  \prod_{i=1}^{n}\left\{ \dfrac{\theta(br^{-1}t^{2-2i} q^{2\mu_i})}
    {\theta(br^{-1}t^{2-2i})}  \left(qt^{2i-2} \right)^{\mu_i} \right\} \\
\cdot \prod_{1\leq i< j \leq n} \hspace*{-5pt} \left\{ \dfrac{
    (qt^{j-i})_{\mu_i - \mu_j} } { (qt^{j-i-1})_{\mu_i - \mu_j} }
    \dfrac{ (br^{-1}t^{3-i-j})_{\mu_i + \mu_j} } {
    (br^{-1}t^{2-i-j})_{\mu_i + \mu_j} } \right\} \\
\cdot W_{\mu} (q^{\lambda}t^{\delta(n)}; q, p, t, bt^{2-2n}, br^{-1}t^{1-n})
\end{multline}
where $x,r,q,p,t, a,b\in\mathbb{C}$.

Note that $W_{\lambda/\mu}(x; q,p, t, a,b)$ vanishes unless $\lambda/\mu$
is a horizontal strip, whereas $\omega_{\lambda/\mu}(x; r;
a,b)=\omega_{\lambda/\mu}(x; r, q,p,t; a,b)$ is defined even when
$\lambda/\mu$ is not a horizontal strip. 

The operator characterization~\cite{Coskun1} of $\omega_{\lambda/\mu}$
yields a recursion formula for \Jc\ in the form 
\begin{equation}
\label{recurrence22}
  \omega_{\lambda/\tau}(y,z; r; a,b) := \sum_\mu
  \omega_{\lambda/\mu}(r^{-k}y; r; ar^{2k},
  b r^k ) \, \omega_{\mu/\tau}(z; r; a, b)
\end{equation}
where $y=(x_{1},\ldots, x_{n-k})\in\mathbb{C}^{n-k}$ and
$z=(x_{n-k+1},\ldots, x_n)\in\mathbb{C}^k$.

A key result used in the development of the multiple special numbers, the \ci\ for
$\omega_{\lambda/\mu}$, is written in~\cite{Coskun1} in the form  
\begin{multline}
\label{cocycleIdentity}
  \omega_{\nu/\mu}((uv)^{-1};uv,q,p,t;a(uv)^2, buv) \\ = \sum_{\mu\subseteq
  \lambda \subseteq \nu} \omega_{\nu/\lambda}(v^{-1};v,q,p,t;a(vu)^2,bvu) \,
  \omega_{\lambda/\mu}(u^{-1};u,q,p,t;au^2,bu)   
\end{multline}
where the summation index $\lambda$ runs over partitions.

Using the recurrence relation~(\ref{recurrence22}) the
definition of $\omega_{\lambda/\mu}(x;r;a,b)$ can be extended 
from the single variable $x\in\mathbb{C}$ case to the multivariable
function $\omega_{\lambda/\mu}(z; r; a,b)$ with arbitrary number of variables
$z = (x_1,\ldots, x_n)\in\mathbb{C}^n$. That $\omega_{\lambda/\mu}(z;
r; a,b)$ is symmetric is also proved in~\cite{Coskun1} using a
remarkable elliptic $BC_n$ \tns\ transformation identity.

Let the $\mathbb{Z}$-space $V$ denote the space of
infinite lower--triangular matrices whose entries are rational
functions 
over the field
$\mathbb{F}=\mathbb{C}(q, p, t,r,a,b)$ as in~\cite{Coskun1}. 
The condition that a matrix $u\in
V$ is lower triangular \wrt\ the partial inclusion
ordering $\subseteq$ defined by
\begin{equation}
\label{partialordering}
\mu \subseteq \lambda \;\Leftrightarrow \;\mu_i \leq \lambda_i, \quad
\forall i\geq 1.
\end{equation}
can be stated in the form
\begin{equation}
  u_{\lambda\mu} = 0 ,\, \quad \mathrm{when}\; \mu \not \subseteq
  \lambda.
\end{equation}
The multiplication operation in $V$ is defined by the relation
\begin{equation}
\label{multiplication}
  (uv)_{\lambda\mu} := \sum_{\mu\subseteq\nu\subseteq\lambda}
  u_{\lambda\nu} v_{\nu\mu}
\end{equation}
for $u,v\in V$. 

\subsection{Limiting Cases}

The limiting cases of the basic (the $p=0$ case of the elliptic) $W$ functions
$W_{\lambda/\mu}(x; q,t,a,b) = W_{\lambda/\mu}(x; q, 0,t,a,b)$ will be
used in computations in what follows. To simplify the exposition, some more notation will be helpful. We set 
\begin{equation}
  W^{ab}_{\lambda/\mu}(x; q,t, s) :=
  \lim_{a\rightarrow 0} W_{\lambda/\mu}(x; q,t, a, as) 
\end{equation}
and 
\begin{equation}
  W^{s\uparrow}_{\lambda/\mu}(x; q,t) :=
  \lim_{s\rightarrow \infty} s^{|\lambda|-|\mu|} W^{ab}_{\lambda/\mu}(x; q,t, s)
\end{equation}
and finally,
\begin{equation}
  W^{s\downarrow}_{\lambda/\mu}(x; q,t) :=
  \lim_{s\rightarrow 0} W^{ab}_{\lambda/\mu}(x; q,t, s)
\end{equation}
The existence of these limits can be seen from ($p=0$ case of) the
definition~(\ref{definitionSkewW}),  the recursion
formula~(\ref{eqWrecurrence}) and the limit rule
\begin{equation}
\label{LimitRule}
  \lim_{a\rightarrow 0}\, a^{|\mu|} (x/a)_{\mu} 
= (-1)^{|\mu|}\, x^{|\mu|} t^{-n(\mu)} q^{n(\mu')} 
\end{equation}
where $\abs{\mu}=\sum_{i=1}^n \mu_i$ and $n(\mu) =
\sum_{i=1}^n (i-1) \mu_i$, 
and $n(\mu') =\sum_{i=1}^n \binom{\mu_i}{2}$.

These functions are closely related to the Macdonald polynomials~\cite{Macdonald1, Okounkov1} and $BC_n$ abelian functions~\cite{Rains1}. 

We now make these definitions more precise. Let $H_{\lambda/\mu}(q,t,b)=H_{\lambda/\mu}(q,0,t,b)$, and define
\begin{multline}
%\label{definitionHfactor}
H_{\lambda/\mu}(q,t) 
:= \lim_{b\rightarrow 0} H_{\lambda/\mu}(q,t,b) \\ 
= \prod_{1\leq i < j\leq
n}\left\{\dfrac{(q^{\mu_i-\mu_{j-1}}t^{j-i})_{\mu_{j-1}-\lambda_j} }
{(q^{\mu_i-\mu_{j-1}+1}t^{j-i-1})_{\mu_{j-1}-\lambda_j} } 
\dfrac{(q^{\lambda_i-\mu_{j-1}+1}t^{j-i-1})_{\mu_{j-1}-\lambda_j}}
{(q^{\lambda_i-\mu_{j-1}}t^{j-i})_{\mu_{j-1}-\lambda_j}}
\right\}
\end{multline}
By setting $b=as$ in the definition of $W_{\lambda/\mu}$ function, and sending $a\rightarrow
0$ we define the family of symmetric rational functions $W^{ab}_{\lambda/\mu}$ in the form
%\begin{equation}\label{definitionSkewW}
%W^{ab}_{\lambda/\mu}(x; q, t, s)
%:= \prod_{i=1}^n t^{i(\mu_i-\lambda_{i+1})} \cdot H_{\lambda/\mu}(q, t) %\dfrac{(x^{-1})_\lambda (qs/(xt))_\mu}
%{(x^{-1})_\mu (qbx, qs/x)_\lambda}
%\end{equation}
% and since \[\prod_{i=1}^n t^{i(\mu_i-\lambda_{i+1})} = t^{-n(\lambda)+|\mu|+n(\mu) } \]
\begin{multline}
W^{ab}_{\lambda/\mu}(x; q, t, s)
:=\lim_{a\rightarrow 0} W_{\lambda/\mu}(x;q,t, a, as) \\ 
= t^{-n(\lambda)+|\mu|+n(\mu) } H_{\lambda/\mu}(q, t) \dfrac{(x^{-1})_\lambda (qs/(xt))_\mu}
{(x^{-1})_\mu (qs/x)_\lambda} \hspace*{59pt}
\end{multline}
for $x\in \mathbb{C}$. Using~(\ref{eqWrecurrence}) 
we get the following recurrence formula for $W^{ab}_{\lambda/\mu}$ function
\begin{equation}
W^{ab}_{\lambda/\mu}(y,z;q, t, s) \\
= \sum_{\nu\prec \lambda} 
W^{ab}_{\lambda/\nu}(yt^{-\ell};q, t, st^{-\ell}) \, 
W^{ab}_{\nu/\mu}(z;q, t, s)
\end{equation}
where $y\in\mathbb{C}$ and $z\in\mathbb{C}^\ell$ as before. 

Similarly, the $W^{s\uparrow}_{\lambda/\mu}$ and $W^{s\downarrow}_{\lambda/\mu}$ are defined as follows. 
\begin{multline}
  W^{s\uparrow}_{\lambda/\mu}(x; q,t) :=
  \lim_{s\rightarrow \infty} s^{|\lambda|-|\mu|} W^{ab}_{\lambda/\mu}(x; q,t, s) \\
  =(-q/x)^{-|\lambda|+|\mu|} q^{ -n(\lambda') + n(\mu') }  H_{\lambda/\mu}(q, t) \dfrac{(x^{-1})_\lambda }{(x^{-1})_\mu } \hspace*{51pt}
\end{multline}
The recurrence formula for $W^{s\uparrow}_{\lambda/\mu}$ function turns out to be
\begin{equation}
\label{Wsuprec}
W^{s\uparrow}_{\lambda/\mu}(y,z;q, t) \\
= \sum_{\nu\prec \lambda} t^{\ell(|\lambda|-|\nu|)}
W^{s\uparrow}_{\lambda/\nu}(yt^{-\ell};q, t) \, 
W^{s\uparrow}_{\nu/\mu}(z;q, t)
\end{equation}
for $y\in\mathbb{C}$ and $z\in\mathbb{C}^\ell$. In the same way, we define
\begin{multline}
  W^{s\downarrow}_{\lambda/\mu}(x; q,t) :=
  \lim_{s\rightarrow 0} W^{ab}_{\lambda/\mu}(x; q,t, s) \\
= t^{-n(\lambda)+|\mu|+n(\mu) } H_{\lambda/\mu}(q, t) \dfrac{(x^{-1})_\lambda }
{(x^{-1})_\mu } \hspace*{107pt}
\end{multline}
The recurrence formula for $W^{s\downarrow}_{\lambda/\mu}$ function becomes
\begin{equation}
W^{s\downarrow}_{\lambda/\mu}(y,z;q, t) \\
= \sum_{\nu\prec \lambda}
W^{s\downarrow}_{\lambda/\nu}(yt^{-\ell};q, t) \, 
W^{s\downarrow}_{\nu/\mu}(z;q, t)
\end{equation}
where again $y\in\mathbb{C}$ and $z\in\mathbb{C}^\ell$. 

\section{$qt$-Binomial Coefficients}
\label{section1}
A common property for all types of special numbers we develop in this paper is that they are  closely connected with binomial coefficients. Therefore, we start with the definition of $qt$-binomial coefficients which will be proved using a multiple analogue of the terminating version of $qt$-binomial theorem. We first derive a multiple terminating $_2\phi_1$ sum as a limit of Jackson's $_8\phi_7$ summation formula. 

\begin{theorem}
For an $n$-part partition $\lambda$, we have
\begin{multline}
\label{2phi1}
 \dfrac{(sx^{-1} )_\lambda } {(s)_\lambda }
 = \sum_{\mu \subseteq \lambda} q^{|\mu|} t^{2n(\mu)} \dfrac{ (x^{-1} )_\mu  } {
   (qt^{n-1} )_\mu} \, 
 \prod_{1\leq i<j \leq 
 n} \left\{\dfrac{ (qt^{j-i})_{\mu_i-\mu_j} }
 {(qt^{j-i-1})_{\mu_i-\mu_j}  }
 \right\}  \\ \cdot W^{ab}_\mu(q^\lambda t^{\delta(n)}; q, t, s^{-1} t^{n-1} ) 
\end{multline}
where $s, x, t, q\in\mathbb{C}$.
\end{theorem}

\begin{proof}
The one variable basic (i.e., $p=0$) version of the $\omega$-\Js~\cite{Coskun1} 
\begin{equation}
  \omega_{\lambda}(s^{-1}x; s; as^2,bs) = \sum_{\mu \subseteq \lambda}
  \omega_{\lambda/\mu}(s^{-1}; s; as^{2},
  b s) \, \omega_{\mu}(x; s; a, b)
\end{equation}
can be written explicitly in the form  
\begin{multline}
\label{simplifiedJs2}
  \dfrac{(sx^{-1}, as x)_\lambda}{(qb x, q b /ax )_\lambda}
 =\sum_{\mu \subseteq \lambda} \dfrac{(s, a
   s)_\lambda} {(qb, qb /a)_\lambda} \dfrac{ (bt^{1-n}, qb/a s)_\mu }{
   (qt^{n-1}, as)_\mu} \\ \cdot\prod_{i=1}^n \left\{\dfrac{(1-b
   t^{2-2i}q^{2\mu_i})} {(1-b t^{2-2i})} (qt^{2i-2})^{\mu_i}\right\}
 \prod_{1\leq i<j \leq 
 n} \left\{\dfrac{ (qt^{j-i})_{\mu_i-\mu_j} (bt^{3-i-j})_{\mu_i+\mu_j}}
 {(qt^{j-i-1})_{\mu_i-\mu_j} (b t^{2-i-j})_{\mu_i+\mu_j} }
 \right\}  \\ \cdot W_\mu(q^\lambda t^{\delta(n)}; q, t, b st^{2-2n},
 bt^{1-n}) \dfrac{(x^{-1}, ax)_\mu} {(qbx, q b /ax)_\mu} 
\end{multline}
Substituting $b=d a$ and sending $a,d\rightarrow 0$ gives the desired result.
\end{proof}
The terminating $qt$-binomial theorem follows as a corollary. 
\begin{cor}
For an $n$-part partition $\lambda$, we have
\begin{multline}
\label{qt_binom_thm}
 (x)_\lambda 
 = \sum_{\mu \subseteq \lambda} (-1)^{|\mu|} q^{|\mu|+n(\mu')} t^{n(\mu)+(1-n)|\mu| } \dfrac{1 }  { (qt^{n-1} )_\mu} \\ \cdot \prod_{1\leq i<j \leq n} \left\{\dfrac{ (qt^{j-i})_{\mu_i-\mu_j} } {(qt^{j-i-1})_{\mu_i-\mu_j}  } \right\}  
  W^{s\uparrow}_\mu(q^\lambda t^{\delta(n)}; q, t) \cdot x^{|\mu|}   
\end{multline}
where $x, t, q\in\mathbb{C}$.
\end{cor}

\begin{proof}
Replace $x\rightarrow sx^{-1}$ and $s\rightarrow sx$ in identity~(\ref{2phi1}) in that order, and send $s\rightarrow 0$ using~(\ref{LimitRule}) to get the identity to be proved.
\end{proof}

This latter identity is a multiple analogue of Cauchy's $q$-binomial theorem~\cite{GasperR1}. Using~(\ref{qt_binom_thm}) we give the definition of a multiple analogue of the binomial coefficient as promised.

\begin{definition}
\label{qtbinomcoeff}
Let $\lambda$ and $\mu$ be $n$-part partitions. Then the $qt$-binomial coefficient is defined by 
\begin{equation}
\binom{\lambda}{\mu}_{\!\!\!q,t} := \dfrac{ q^{|\mu|} t^{2n(\mu)+(1-n)|\mu| } }  { (qt^{n-1} )_\mu} \prod_{1\leq i<j \leq n} \left\{\dfrac{ (qt^{j-i})_{\mu_i-\mu_j} } {(qt^{j-i-1})_{\mu_i-\mu_j}  } \right\}  W^{s\uparrow}_\mu(q^\lambda t^{\delta(n)}; q, t)
\end{equation}
where $q,t\in\mathbb{C}$.
\end{definition}
Note that with this definition we can write the terminating $qt$-binomial theorem~(\ref{qt_binom_thm}) in the form
\begin{equation}
\label{qt_binom_thmAlt}
 (x)_\lambda 
 = \sum_{\mu \subseteq \lambda} (-1)^{|\mu|} q^{n(\mu')} t^{-n(\mu)} \binom{\lambda}{\mu}_{\!\!\!q,t} x^{|\mu|} 
\end{equation}
Note also that setting $t=q^\alpha$ and sending $q\rightarrow 1$ yields multiple ordinary $\alpha$-binomial coefficients. Below we extend this definition further to be valid not only for partitions $\lambda$ and $\mu$, but also for any $n$-tuples of complex numbers  $\lambda\in\mathbb{C}^n$ and $\mu\in\mathbb{C}^n$. For $n=1$, the definition reduces to that of the one dimensional $q$-binomial coefficients
\begin{equation}
\binom{n}{k}_{\!\!q} : = \dfrac{ (q)_n }{(q)_{n-k} (q)_k }
\label{eq:qBinom}
\end{equation}
which are also known as the Gaussian polynomials that are studied extensively in the literature including but not limited to the works in~\cite{Andrews5, Andrews6, GasperR1, BerkovichW1, Gould2, Macdonald2, GarvanS1, ClarkI1}. 

Some of the main properties of the binomial coefficients readily generalize to the multiple case. For example, the identities
\begin{equation}
2^n=\sum_{k=0}^n \binom{n}{k} \quad\mathrm{and}\quad
0=\sum_{k=0}^n (-1)^k \binom{n}{k}
\end{equation}
have the following multiple analogues.

\begin{theorem}
For an $n$-part partition $\lambda$
\begin{equation}
 (-1)_\lambda 
 = \sum_{\mu \subseteq \lambda} q^{n(\mu')} t^{-n(\mu)} \binom{\lambda}{\mu}_{\!\!\!q,t}  
\end{equation}
and
\begin{equation}
 0=  \sum_{\mu \subseteq \lambda} (-1)^{|\mu|} q^{n(\mu')} t^{-n(\mu)} \binom{\lambda}{\mu}_{\!\!\!q,t}
\end{equation}
\end{theorem}
\begin{proof}
These identities follow immediately from~(\ref{qt_binom_thmAlt}) by setting $z=-1$ and $z=0$, respectively.
\end{proof}

Before investigating multiple analogues of other binomial identities, we first note that the definition~(\ref{qtbinomcoeff}) makes sense even for generalized partitions $\mu=(\mu_1,\mu_2,\ldots, \mu_n)$ such that $\mu_1\geq \mu_2\geq \cdots \geq \mu_n$ with possibly negative parts $\mu_i$ starting with some $i\in[n]:=\{1,2,\ldots,n \}$. 

\begin{theorem}
\label{thm:extension}
The $W_\lambda$ function is well-defined defined for the generalized partitions $\lambda$. 
In fact, the evaluation of $W_\lambda$ function that enters into the definition of $qt$-binomial coefficient can be computed for $n$-tuples of complex numbers $\lambda\in\mathbb{C}^n$ and $\mu\in\mathbb{C}^n$. The limiting $W^{ab}_\lambda$, $W^{s\uparrow}_\lambda$ and $W^{s\downarrow}_\lambda$ functions are also defined in the general case. 
\end{theorem}
\begin{proof}
Let $\lambda$ be an $n$-part partition with $\lambda_n 
\neq 0$ and $0\leq k\leq \lambda_n$ for some integer $k$, and let 
$z=(z_1,\ldots,z_n)\in \mathbb{C}^n$. We have shown~\cite{CoskunG1} that
\begin{multline}
\label{commonfactor}
W_{\lambda}(z;q,t,a,b) \\
=\prod_{1\leq i<j\leq n}\frac{(qbt^{j-2i})_{2k}}{(qbt^{j-1-2i})_{2k}}
\prod_{i=1}^n\frac{(z_i^{-1})_k (az_i)_k}{(qbz_i)_k (qb/(az_i))_k}
W_{\lambda-k^n}(zq^{-k}; aq^{2k}, bq^{2k}, t, q)
\end{multline}
Among other applications, this identity extends the definition of $W_{\lambda}$ function to general partitions. For a generalized partition $\lambda$, the index $\lambda-k^n$ on the \rhs\ will be a standard partition.

The duality formula for $W_{\lambda}$ functions from~\cite{CoskunG1} may be used to compute the special evaluation of $W_{\lambda}$ that occur in many application in this paper even for $\lambda\in\mathbb{C}^n$. The formula can be stated as 
\begin{multline}
\label{eq:duality}
W_{\lambda}\left(k^{-1}q^\nu t^\delta;q,t,k^2a,kb\right)
\cdot\dfrac{(qbt^{n-1})_\lambda (qb/a)_\lambda}{(k)_\lambda
  (kat^{n-1})_\lambda} \\
\cdot \prod_{1\leq i < j\leq n} \left\{ \dfrac{(t^{j-i})_{\lambda_i
    -\lambda_j}(qa^{\prime}t^{2n-i-j-1})_{\lambda_i+\lambda_j}}
{(t^{j-i+1})_{\lambda_i
-\lambda_j}(qa^{\prime}t^{2n-i-j})_{\lambda_i+\lambda_j}} \right\}\\
=W_{\nu}\left(h^{-1}q^\lambda t^\delta;q, t,h^2a^{\prime},hb\right)
\cdot \dfrac{(qbt^{n-1})_\nu (qb/a^{\prime})_\nu}{(h)_\nu
  (ha^{\prime}t^{n-1})_\nu} \\
\cdot \prod_{1\leq i < j\leq n} \left\{ \dfrac{(t^{j-i})_{\nu_i
      -\nu_j}(qat^{2n-i-j-1})_{\nu_i+\nu_j}} {(t^{j-i+1})_{\nu_i
-\nu_j}(qat^{2n-i-j})_{\nu_i+\nu_j}} \right\}
\end{multline}
where $k=a^{\prime}t^{n-1}/b$ and $h = at^{n-1}/b$. 
Substitute these latter relations between parameters into the identity 
and factor out one term
\begin{equation}
\dfrac{q^{-\lambda_1\nu_1} (k)_{\lambda_1} (q/k )_{\nu_1} (k a t^{2(n-1)} q^{\nu_1}  )_{\lambda_1} }
{(q^{1-\lambda_1}/k )_{\nu_1}  (b t^{n-1} q^{1+\nu_1} , (b t^{1-n}/a) q^{1-\nu_1} )_{\lambda_1} } 
\end{equation}
corresponding to the dominant weight $\lambda_1$ from the $W$ function on the \lhs, and one term 
\begin{equation}
\dfrac{((at^{n-1}/b) q^{-\lambda_1}, a k t^{2(n-1)} q^{\lambda_1} )_{\nu_1} }
{( bt^{n-1}  q^{1+\lambda_1}, q^{1-\lambda_1}/k ))_{\nu_1} } 
\end{equation}
for $\nu_1$ from the $W$ on the \rhs\ to get
\begin{multline}
\dfrac{q^{-\lambda_1\nu_1} (k)_{\lambda_1} (q/k )_{\nu_1} (k a t^{2(n-1)} q^{\nu_1}  )_{\lambda_1} }
{(q^{1-\lambda_1}/k )_{\nu_1}  (b t^{n-1} q^{1+\nu_1} , (b t^{1-n}/a) q^{1-\nu_1} )_{\lambda_1} } \, 
W'_{\lambda}\left(k^{-1}q^\nu t^\delta;q, t, k^2a, kb\right) \\
\cdot\dfrac{(qbt^{n-1})_\lambda (qb/a)_\lambda}{(k)_{\lambda_1} (kt^{-1})_{\lambda'}
  (kat^{n-1})_\lambda} 
\prod_{1\leq i < j\leq n} \left\{ \dfrac{(t^{j-i})_{\lambda_i
    -\lambda_j}(qk b t^{n-i-j})_{\lambda_i+\lambda_j}}
{(t^{j-i+1})_{\lambda_i
-\lambda_j}(qk b t^{n-i-j+1})_{\lambda_i+\lambda_j}} \right\}\\
=\dfrac{((at^{n-1}/b) q^{-\lambda_1}, a k t^{2(n-1)} q^{\lambda_1} )_{\nu_1} }
{(q^{1-\lambda_1}/k )_{\nu_1} ( bt^{n-1}  q^{1+\lambda_1})_{\nu_1}  } 
\, W'_{\nu}\left((bt^{1-n}/a) q^\lambda t^\delta; q, t, 
a^2 k t^{n-1}/b , at^{n-1}\right) \\
\cdot \dfrac{(qbt^{n-1})_\nu (q t^{n-1}/k )_\nu}{(at^{n-1}/b)_\nu
  ( a k t^{n-1} )_\nu} 
\prod_{1\leq i < j\leq n} \left\{ \dfrac{(t^{j-i})_{\nu_i
      -\nu_j}(qat^{2n-i-j-1})_{\nu_i+\nu_j}} {(t^{j-i+1})_{\nu_i
-\nu_j}(qat^{2n-i-j})_{\nu_i+\nu_j}} \right\}
\end{multline}
where $W'$ denotes the $W$ function after the dominant factors are taken out and $\lambda'$ denotes the partition $(\lambda_2, \lambda_3, \ldots, \lambda_n )$. Setting $k=1$, moving entries to the \rhs, and multiplying back both sides by the dominant factors for the \lhs\ gives
\begin{multline}
W_{\lambda}\left(q^\nu t^\delta;q, t, a, b\right) \\
= \dfrac{q^{\lambda_1\nu_1}  (q^{-\nu_1} )_{\lambda_1} }
{(q )_{\nu_1}  } \dfrac{(t^{-1})_{\lambda'}
  (at^{n-1})_\lambda }{ (qbt^{n-1})_\lambda (qb/a)_\lambda } 
\prod_{1\leq i < j\leq n} \left\{ \dfrac{(t^{j-i+1})_{\lambda_i
-\lambda_j}(q b t^{n-i-j+1})_{\lambda_i+\lambda_j}}{(t^{j-i})_{\lambda_i
    -\lambda_j}(q b t^{n-i-j})_{\lambda_i+\lambda_j}}
 \right\} \\
\cdot \dfrac{((at^{n-1}/b) q^{-\lambda_1}, a t^{2(n-1)} q^{\lambda_1} )_{\nu_1} }
{ ( bt^{n-1}  q^{1+\lambda_1})_{\nu_1}  } 
\, W'_{\nu}\left((bt^{1-n}/a) q^\lambda t^{\delta(n)}; q, t, 
a^2 t^{n-1}/b , at^{n-1}\right) \\
\cdot \dfrac{(qbt^{n-1})_\nu (q t^{n-1} )_\nu}{(at^{n-1}/b)_\nu
  ( a t^{n-1} )_\nu} 
\prod_{1\leq i < j\leq n} \left\{ \dfrac{(t^{j-i})_{\nu_i
      -\nu_j}(qat^{2n-i-j-1})_{\nu_i+\nu_j}} {(t^{j-i+1})_{\nu_i
-\nu_j}(qat^{2n-i-j})_{\nu_i+\nu_j}} \right\}
\end{multline}
This formula allows us to extend the evaluation of the $W_\lambda$ function on the \lhs\ from an $n$-part partition $\lambda$ to any $\lambda\in\mathbb{C}^n$, because $\lambda$ appears as a variable on the \rhs. The front factors that involve $\lambda$ can be extended to the complex case by the definition of $q$-Pochhammer symbol~(\ref{qPochSymbol}).

Similar evaluations for the limiting cases $W^{ab}_\lambda$, $W^{s\uparrow}_\lambda$ and $W^{s\downarrow}_\lambda$ can easily be computed from this result. 
\end{proof}

We now turn to the study of main properties of $qt$-binomial coefficients. 
\begin{theorem}
\label{initial}
For an $n$-part partition $\lambda$ and generalized $n$-part partition $\mu$, we have
\begin{equation}
\binom{\lambda}{\lambda}_{\!\!\!q,t} = \binom{\lambda}{0^n}_{\!\!\!q,t} = 1. 
\end{equation}
where $0^n$ is the $n$-part partitions whose parts are all zero. In addition, 
\begin{equation}
\binom{\lambda}{\mu}_{\!\!\!q,t} = 0 
\end{equation}
when $\mu\not\subseteq \lambda$ or $0^n\not\subseteq\mu$.
\end{theorem}
\begin{proof}
Recall that the normalizing coefficients for $W_\lambda$ function is written~\cite{CoskunG1} in the form
\begin{multline}
\label{Wnormal}
W_\lambda(q^\lambda t^{\delta(n)}; q,
t, a, b) \\
= \prod_{k=1}^n\left\{\frac{(qbt^{n-k}, qt^{n-k})_{\lambda_k}
    (at^{2n-2k})_{2\lambda_k}}
{((a/b)t^{n-k}, at^{n-k})_{\lambda_k} (qbt^{n+1-2k})_{2\lambda_k}}
t^{(n+1-2k)\lambda_k}\right\}\\
\cdot(a/(qb))^{|\lambda|}\cdot\prod_{1\leq i < j\leq
n}\frac{(qt^{j-i-1})_{\lambda_i-\lambda_j}(at^{2n-i-j})_{\lambda_i+\lambda_j}}
{(qt^{j-i})_{\lambda_i-\lambda_j}(at^{1+2n-i-j})_{\lambda_i+\lambda_j}}
\end{multline}
Setting $b=as$ and sending $a\rightarrow 0$ 
and $s\rightarrow\infty$ (after multiplying by $s^{|\lambda|}$) gives
\begin{equation}
\label{Wsnormal}
W^{s\uparrow}_\lambda(q^\lambda t^{\delta(n)}; q, t) 
= ( qt^{n-1} )_{\lambda} \, t^{(n-1)|\lambda|-2n(\lambda)} 
q^{-|\lambda|}\, \prod_{1\leq i < j\leq
n}\frac{(qt^{j-i-1})_{\lambda_i-\lambda_j} }
{(qt^{j-i})_{\lambda_i-\lambda_j} }
\end{equation}
Substituting this evaluation into $\binom{\lambda}{\lambda}_{\!\!q,t}$ gives 1 as desired. 

It is easy to see from the definition of $W_\lambda$ function that $W_\mu(x; q, t, a, b)=1$ when $\mu=0^n$. Therefore $W^s_{0^n}(x; q, t)=1$ and $\binom{\lambda}{\mu}_{\!\!q,t}=1$.

Finally, the fundamental vanishing property of $W$ functions states~\cite{CoskunG1} that 
\begin{equation}
\label{FundamentalVanishing}
  W_{\mu}(q^\lambda t^\delta;q, t,a,b)=0
\end{equation}
when $\mu\not\subseteq \lambda$. Hence $W^{s}_\mu(q^\lambda t^{\delta(n)}; q, t)=0$ and therefore $\binom{\lambda}{\mu}_{\!\!q,t}\!\!\!=0$ in that case. Note also that if $\mu_n<0$ then  $\binom{\lambda}{\mu}_{\!\!q,t}\!\!\! = 0$ due to the $ (qt^{n-1} )_\mu $ factor in the denominator.
\end{proof}

We introduce some notation before writing a multiple analogue of the important recurrence relation 
\[\binom{n+1}{k} = \binom{n}{k} + \binom{n}{k-1}, \quad 1\leq k\leq n \]
for classical binomial coefficients. Let $e_i=(0,\ldots, 1, \ldots,0)$ be the $n$-dimensional vector whose $i$-th coordinate is $1$ and all others are 0. Let $\lambda^i=\lambda+e_i$ whenever $\lambda^i$ is a partition. Finally, we write $\lambda \vdash k$ to denote that $\lambda$ is a partition of the positive integer $k$.

\begin{theorem}
\label{binomialRec}
For an $n$-part partition $\lambda$ and $k\leq |\lambda|$, we have
\begin{multline}
\sum_{\mu\vdash k} q^{n(\mu')} t^{-n(\mu)} \binom{\lambda^i}{\mu}_{\!\!\!q,t} \\
 =
 \sum_{\tau \vdash k} q^{n(\tau')} t^{-n(\tau)} \binom{\lambda}{\tau}_{\!\!\!q,t} + 
\sum_{\nu \vdash (\!k-1\!)} \!\!\! q^{n(\nu')+\lambda_i } t^{-n(\nu)+1-i} \binom{\lambda}{\nu}_{\!\!\!q,t}
\end{multline}
where $\mu\subseteq\lambda^i$ and $\tau,\nu\subseteq\lambda$. 
\end{theorem}
\begin{proof}
The proof is a simple application of the $qt$-binomial theorem~(\ref{qt_binom_thmAlt}). Since
$ (x)_{\lambda^i} = (1- x t^{1-i} q^{\lambda_i}) (x)_{\lambda} $, we get
\begin{multline}
 \sum_{\mu \subseteq \lambda^i} (-1)^{|\mu|} q^{n(\mu')} t^{-n(\mu)} \binom{\lambda^i}{\mu}_{\!\!\!q,t} x^{|\mu|}  \\
 = (1- x t^{1-i} q^{\lambda_i})  \sum_{\tau \subseteq \lambda} (-1)^{|\tau|} q^{n(\tau')} t^{-n(\tau)} \binom{\lambda}{\tau}_{\!\!\!q,t} x^{|\tau|} 
\end{multline}
by a double application of~(\ref{qt_binom_thmAlt}). 
This may be written as
\begin{multline}
 \sum_{k=0}^{|\lambda^i|} \sum_{\mu \vdash k} (-1)^{|\mu|} q^{n(\mu')} t^{-n(\mu)} \binom{\lambda^i}{\mu}_{\!\!\!q,t} x^{|\mu|}  \\
 = (1- x t^{1-i} q^{\lambda_i})  \sum_{\ell =0}^{|\lambda|} \sum_{\tau \vdash \ell} (-1)^{|\tau|} q^{n(\tau')} t^{-n(\tau)} \binom{\lambda}{\tau}_{\!\!\!q,t} x^{|\tau|} 
\end{multline}
Note that the coefficient of $x^k$ on the \lhs\ is 
\begin{equation}
\sum_{\mu \vdash k} (-1)^{|\mu|} q^{n(\mu')} t^{-n(\mu)} \binom{\lambda^i}{\mu}_{\!\!\!q,t} 
\end{equation}
The coefficient of $x^k$ that comes from the two pieces on the \rhs\ becomes 
\begin{equation}
\sum_{\tau \vdash k} (-1)^{|\tau|} q^{n(\tau')} t^{-n(\tau)} \binom{\lambda}{\tau}_{\!\!\!q,t} 
 +\!\! \sum_{\nu \vdash (\!k-1\!)} (-1)^{|\nu|+1} q^{n(\nu')+\lambda_i} t^{-n(\nu)+1-i} \binom{\lambda}{\nu}_{\!\!\!q,t}
\end{equation}
Note also that $\mu\subseteq\lambda^i$ and $\tau,\nu\subseteq\lambda$, for otherwise the binomial coefficients vanish by Theorem~(\ref{initial}) above. Canceling out common factors gives the result. 
\end{proof}

A multiple analogue of the symmetry property 
\[\binom{n}{k} = \binom{n}{n-k}, \quad 0\leq k\leq n \]
for classical binomial coefficients now follows.
\begin{cor}
For an $n$-part partition $\lambda$ and $k\leq |\lambda|$, we have
\begin{equation}
\sum_{\mu\vdash k} \binom{\lambda}{\mu}_{\!\!\!q,t} = \sum_{\tau \vdash (|\lambda|-k)} \binom{\lambda}{\tau}_{\!\!\!q,t}
\end{equation}
\end{cor}
\begin{proof}
This result follows from the $qt$-binomial theorem~(\ref{qt_binom_thmAlt}) through a similar argument used in the proof of the previous theorem.
\end{proof}

We will now write an analogue of the identity
\begin{equation}
\label{inv}
\dfrac{z^k}{(1-z)^{k+1} } = \sum_{n=k}^\infty
  \binom{n}{k} z^n,  \quad |z|<1
\end{equation}
Let $\lambda$ be an $n$-part partition and consider the $\mathbb{C}(q,t)$-space $P_\lambda$ of all polynomials of degree less then or equal to $|\lambda|$ in a single complex variable $x$. Note that both sets $\beta_1=\{x^k: k=1,\ldots, |\lambda| \}$ and $\beta_2=\{(x)_{\mu}: \mu\subseteq \lambda \}$ form bases for $P_\lambda$. In fact, $qt$-binomial theorem stipulates just this, providing a change of basis formula between the two bases. 

Similarly, we view the terminating $_2\phi_1$ sum~(\ref{2phi1}) as a matrix representation of the shift operator acting on $\beta_2$, and prove an simpler version of the \ci~\cite{CoskunG1} for $\omega_\lambda$ functions as follows.
\begin{theorem}
For $n$-part partitions $\nu$ and $\mu$, we have
\begin{multline}
\label{weakci}
(sr)_\nu
 \, W^{ab}_\mu(q^\nu t^{\delta(n)}; q, t, (sr)^{-1} t^{n-1} ) \\
 = \sum_{\mu \subseteq \lambda \subseteq \nu} 
 q^{|\lambda|} t^{2n(\lambda)} 
 \dfrac{(s)_\nu} { (qt^{n-1} )_\lambda} 
 \prod_{1\leq i<j \leq 
 n} \left\{\dfrac{ (qt^{j-i})_{\lambda_i-\lambda_j} }
 {(qt^{j-i-1})_{\lambda_i-\lambda_j}  }
 \right\}  W^{ab}_\lambda(q^\nu t^{\delta(n)}; q, t, s^{-1} t^{n-1} ) \\
 \cdot (r)_\lambda \,
   W^{ab}_\mu(q^\lambda t^{\delta(n)}; q, t, r^{-1} t^{n-1} ) 
\end{multline}
\end{theorem}

\begin{proof}
Write the identity~(\ref{2phi1}) in the form
\begin{multline}
 (sx )_\lambda
 = \sum_{\mu \subseteq \lambda} q^{|\mu|} t^{2n(\mu)} 
 % \prod_{i=1}^n \left\{ (qt^{2i-2})^{\mu_i}\right\} 
 \dfrac{(s)_\lambda} { (qt^{n-1} )_\mu} 
 \prod_{1\leq i<j \leq 
 n} \left\{\dfrac{ (qt^{j-i})_{\mu_i-\mu_j} }
 {(qt^{j-i-1})_{\mu_i-\mu_j}  }
 \right\}  \\ \cdot W^{ab}_\mu(q^\lambda t^{\delta(n)}; q, t, s^{-1} t^{n-1} ) 
  (x)_\mu  
\end{multline}
Starting with the basis $\beta_2$, apply the shift operator by a factor of $r$ followed by a shift by $s$. This double shift can be achieved by a composite shift by a factor of $sr$. Writing this argument explicitly using the identity above and simplifying gives the desired result.
\end{proof}

Using this result, we now write a multiple analogue of the binomial identity~(\ref{inv}) given above. Here $(a)_{\infty^n}$ denotes $\prod_{i=1}^n (a t^{1-i})_\infty $. 
\begin{theorem}
For $n$-part partitions $\nu$ and $\mu$, we have
\begin{equation}
\dfrac{z^{|\mu|}}{ (qz)_{\infty^n} }  
\prod_{1\leq i < j\leq n} \dfrac{(t^{j-i+1})_{\mu_i
-\mu_j} } {(t^{j-i})_{\mu_i -\mu_j} } \\
 = \sum_{\lambda \supseteq  \mu} 
 z^{|\lambda|}
 \prod_{1\leq i<j \leq 
 n} 
 \dfrac{(t^{j-i+1})_{\lambda_i
-\lambda_j} } {(t^{j-i})_{\lambda_i -\lambda_j} }
  \cdot \binom{\lambda}{\mu}_{\!\!\!q,t}
\end{equation}
where $\max \{\abs{ qzt^{(2i-n-1)} } : \, i\in[n] \} < 1$. 
\end{theorem}
 
\begin{proof}
First send $r\rightarrow 0$ in the~(\ref{weakci}) to get 
\begin{multline}
\label{step1}
W^{s}_\mu(q^\nu t^{\delta(n)}; q, t) \\
 = \sum_{\mu \subseteq \lambda \subseteq \nu} 
 q^{|\lambda|} t^{2n(\lambda)} 
 \dfrac{ (s)_\nu } { (qt^{n-1} )_\lambda} 
 \prod_{1\leq i<j \leq 
 n} \left\{\dfrac{ (qt^{j-i})_{\lambda_i-\lambda_j} }
 {(qt^{j-i-1})_{\lambda_i-\lambda_j}  }
 \right\}  \! W^{ab}_\lambda(q^\nu t^{\delta(n)}; q, t, s^{-1} t^{n-1} ) \\
  \cdot  s^{-|\mu|} W^{s}_\mu(q^\lambda t^{\delta(n)}; q, t ) 
\end{multline}
Next we send $\nu\rightarrow \infty$. Note that setting $\nu=k^n$, where $k^n=(k,k,\ldots, k)$ denotes the $n$-part partition whose parts all equal $k$, simplifies the $W$ function evaluations in the identity. To make this precise, we recall the Weyl denominator formula~\cite{CoskunG1} which states that
\begin{multline}
\label{Weyl}
W_{\mu}(xt^{\delta(n)};q,t,a,b) \\
=\dfrac{(x^{-1}, axt^{n-1})_\mu}{(qbxt^{n-1}, qb/(ax))_\mu} 
\! \prod_{1\leq i < j\leq n} \dfrac{(t^{j-i+1})_{\mu_i
-\mu_j}(qbt^{n-i-j+1})_{\mu_i+\mu_j}}
{(t^{j-i})_{\mu_i -\mu_j}(qbt^{n-i-j})_{\mu_i+\mu_j}}
\end{multline}
Set $b=cs^{-1}t^{1-n}$ and $a=ct^{2-2n}$ in~(\ref{Weyl}), and send $c\rightarrow 0$ to get
\begin{equation}
\label{Wablimit}
W^{ab}_{\mu}(xt^{\delta(n)};q,t,s^{-1}t^{n-1}) \\
=\dfrac{(x^{-1})_\mu}{(qs^{-1}t^{n-1}/x)_\mu} 
\! \prod_{1\leq i < j\leq n} \dfrac{(t^{j-i+1})_{\mu_i
-\mu_j} } {(t^{j-i})_{\mu_i -\mu_j} }
\end{equation}
Multiplying both sides by $(s^{-1}t^{n-1})^{|\mu|}$, and sending $s\rightarrow 0$ further implies that
\begin{equation}
\label{Wslimit}
W^{s\uparrow}_{\mu}(xt^{\delta(n)};q,t) \\
=  (-1)^{|\mu|} x^{|\mu|} t^{n(\mu)} q^{-|\mu|-n(\mu')} (x^{-1})_\mu
\! \prod_{1\leq i < j\leq n} \dfrac{(t^{j-i+1})_{\mu_i
-\mu_j} } {(t^{j-i})_{\mu_i -\mu_j} }
\end{equation}
On setting $x=q^k$, sending $k\rightarrow \infty$ 
in the last two evaluations and substituting them into the identity~(\ref{step1}), we get 
\begin{multline}
\dfrac{s^{|\mu|} q^{-|\mu|}}{ (s)_{\infty^n} }  
\prod_{1\leq i < j\leq n} \dfrac{(t^{j-i+1})_{\mu_i
-\mu_j} } {(t^{j-i})_{\mu_i -\mu_j} } 
 = \sum_{\lambda \supseteq \mu } 
 \dfrac{t^{2n(\lambda)-(n-1)|\lambda|}  s^{|\lambda|} } 
 { (qt^{n-1} )_\lambda} \\ \cdot
 \prod_{1\leq i<j \leq 
 n} \left\{\dfrac{ (qt^{j-i})_{\lambda_i-\lambda_j} }
 {(qt^{j-i-1})_{\lambda_i-\lambda_j}  } 
 \dfrac{(t^{j-i+1})_{\lambda_i
-\lambda_j} } {(t^{j-i})_{\lambda_i -\lambda_j} }
 \right\}    W^{s\uparrow}_\mu(q^\lambda t^{\delta(n)}; q, t ) 
\end{multline}
Replacing $s$ by $qz$ gives the identity to be proved upon verifying the convergence.

The fact that the resulting infinite series converges follows from a multiple analogue of the dominated convergence theorem introduced in~\cite{Coskun1}. Consider the multiple series
\begin{equation}
  \sum_{\lambda\in L^+_k} h_\lambda(k)
\end{equation}
where $L^+_k$ is the sub--alcove $\{\lambda\in \mathbb{Z}^n: k\geq
\lambda_1 \geq \ldots \lambda_n \geq 0 \}$. The theorem states that if the pointwise limit $ h_\lambda:= \lim_{k\rightarrow\infty}
h_{\lambda}(k)$ exists for all $\lambda\in L^+_k$, and we can find
$m^h_\lambda$ for each $\lambda$ such that $\abs{h_\lambda(k)}
\leq m^h_\lambda$ for all $k \geq \lambda_1$, and that the series $\sum_{\lambda\in L^+} m^h_\lambda$ is convergent, then the original series converges.

That the pointwise limit exists on both sides in~(\ref{step1}) is already verified above. Note  that the index $\mu$ of the $W$ function inside the summand is a fixed partition. Note also that, except the powers of $q,t$ and $s$, all other factors in the summand of~(\ref{step1}) can be put into the form $(uq^{\alpha})_\infty /(vq^{\alpha})_\infty$ using the definition of $q$-Pochhammer symbol~(\ref{qPochSymbol}). 
Standard theorems on infinite products and sequences imply that such factors
are bounded when $\alpha$ is a non--negative integer,  
$u,v\in\mathbb{C}$, and that $v$ is such that the denominator never
vanishes. This is because
\begin{equation}
  \lim_{\alpha \rightarrow \infty} \left\vert \dfrac{
     (uq^{\alpha})_\infty } {(vq^{\alpha})_\infty } \right\vert = 1
\end{equation}
when $\abs{q}<1$.
Therefore it follows that  
for a constant  $C_h$ that may depend only on $q, t$ and $s$,
and is independent of $k$ and $\lambda$, we've found
\begin{equation}
 m^h_\lambda = C_h \abs{ t^{2n(\lambda)-(n-1)|\lambda|} s^{|\lambda|} } =  C_h \abs{ \prod_{i=1}^n t^{(2i-n-1)\lambda_i} s^{\lambda_i} }     
\end{equation}
Finally, we need to show that $\sum_{\lambda\in L^+} m^h_\lambda$ is
convergent which may be verified using the multiple series ratio
test. Let $e_i$ be defined as above. Then, we see that 
\begin{equation}
 \abs{ \dfrac{ m^h_{\lambda+e_i } } {m^h_\lambda } } = \abs{ t^{(2i-n-1)} qz}     
\end{equation}
for each $i\in [n]$ where possible (i.e., when $\lambda+e_i$ is a partition). Therefore, $\sum_{\lambda\in L^+} m^h_\lambda$ converges when 
\[ \max_{i\in[n]} \{\abs{ t^{(2i-n-1)} qz}  \} < 1.\]
The original series also converges under the same condition.
\end{proof}

Next, we write a multiple analogue of the identity 
\begin{equation}
2^{n-m} \binom{n}{m} = \sum_{k=m}^n 
   \binom{n}{k} \binom{k}{m}
\end{equation}
for $qt$-binomial coefficients. 
\begin{theorem}
For $n$-part partitions $\nu$ and $\mu$, we have
\begin{equation}
\label{doublebinom}
t^{-n(\mu)} q^{n(\mu')} 
\dfrac{ ( -1 )_{\nu}} { (-1 )_{\mu} } \,\binom{\nu}{\mu}_{\!\!\!q,t} 
= \sum_{\mu \subseteq
    \lambda \subseteq \nu } 
    t^{-n(\lambda)} q^{n(\lambda')} 
    \binom{\nu}{\lambda}_{\!\!\!q,t}
    \binom{\lambda}{\mu}_{\!\!\!q,t}
\end{equation}
\end{theorem}
 
\begin{proof}
The \ci~\cite{CoskunG1} for $\omega$ functions
\begin{equation}
  \omega_{\lambda/\mu}((sr)^{-1},sr,as^2, bs) 
=\sum_\nu
  \omega_{\lambda/\nu}(s^{-1},s,as^2,bs) \,
  \omega_{\nu/\mu}(r^{-1};r,a,b)  
\end{equation}
may be written explicitly in the form 
\begin{multline}
\dfrac{(rs)_{\lambda} (asr^{-1})_{\lambda}} {(qbr^{-1})_{\lambda}
    (qbr/a)_{\lambda}} \dfrac {(qb)_{\lambda}
    (qb/a)_{\lambda}} {(s)_{\lambda} (as)_{\lambda}}
    \dfrac{ (qb/as)_{\mu}}{ (asr^{-1})_{\mu}}  \dfrac{
    (ar^{-1})_{\mu}} {(qb/a)_{\mu}}  \\      
 W_{\mu} (t, bst^{2-2n}, br^{-1}t^{1-n}; q^{\lambda}t^{\delta(n)};q) \\
=  \sum_\nu \dfrac{ (qb/as)_{\nu}}{ (as)_{\nu}} \dfrac{
    (bt^{1-n})_{\nu}}{ (qt^{n-1})_{\nu}} \dfrac{(r)_{\nu}
    (ar^{-1})_{\nu}} {(qbr^{-1})_{\nu} (qbr/a)_{\nu}}\\     
  \prod_{i=1}^{n}\left\{ \dfrac{(1-bt^{2-2i} q^{2\nu_i})}
    {(1-bt^{2-2i})}  \left(qt^{2i-2} \right)^{\nu_i} \right\} 
  \prod_{1\leq i< j \leq n} \left\{ \dfrac{ (qt^{j-i})_{\nu_i - \nu_j}
    } { (qt^{j-i-1})_{\nu_i - \nu_j} } \dfrac{
    (bt^{3-i-j})_{\nu_i + \nu_j} } { (bt^{2-i-j})_{\nu_i + 
    \nu_j} } \right\} \\ W_{\nu} (t, bst^{2-2n}, bt^{1-n};
    q^{\lambda}t^{\delta(n)};q) \, W_{\mu} (t, bt^{2-2n}, br^{-1}t^{1-n};
    q^{\nu}t^{\delta(n)};q)   
\end{multline}
Setting $b=azrs/q$, and sending $a\rightarrow 0$ and $r\rightarrow 0$ gives
\begin{multline}
\label{strongci_limit}
\dfrac{ ( zs )_{\nu}} {(s )_{\nu}} 
    \dfrac{( z )_{\mu}}{ (zs )_{\mu} } s^{|\mu|} 
W^{s}_{\mu} (q^{\nu}t^{\delta(n)}; q, t) \\
=  \sum_{\mu \subseteq
\lambda \subseteq \nu} \dfrac{(z )_{\lambda}}{(qt^{n-1})_{\lambda} } 
\, q^{|\lambda|} t^{2n(\lambda)} 
\cdot \prod_{1\leq i< j \leq n} \hspace*{-5pt} \left\{ \dfrac{
    (qt^{j-i})_{\lambda_i - \lambda_j} } { (qt^{j-i-1})_{\lambda_i - \lambda_j} }
      \right\} \\
\cdot W^{ab}_{\lambda} (q^{\nu}t^{\delta(n)}; q, t, s^{-1}t^{n-1} ) 
\cdot W^{s\uparrow}_{\mu} (q^{\lambda}t^{\delta(n)}; q, t) 
\end{multline}
We will use a multiple analogue of Bailey's $_{10}\phi_9$ transformation formula from~\cite{CoskunG1} to transform the series on the \rhs. 
The $_{10}\phi_9$ transformation can be written explicitly in the form
\begin{multline}
\dfrac{(v)_{\nu} (a's^2 v^{-1})_{\nu}} {(qbrs/v)_{\nu}
    (qbrv/sa')_{\nu}} \dfrac {(qbr)_{\nu}
    (qbr/a')_{\nu}} {(s)_{\nu} (a's)_{\nu}} \dfrac{ (Av^2/s^2r)_{\mu}}{
    (qbrs^2/Av^2)_{\mu}}  \dfrac {(qbrs/Av)_{\mu}}{ (Av/sr)_{\mu}} \\ \cdot
  \sum_{\substack{\lambda\\ \mu \subseteq
    \lambda \subseteq \nu }} \dfrac{ (qbrs^2/Av^2)_{\lambda}}{
    (Av/s)_{\lambda}} \dfrac{ 
    ((brsv^{-1}t^{1-n})_{\lambda}}{ (qt^{n-1})_{\lambda}} \dfrac{(rs/v)_{\lambda} 
    (Av/sr )_{\lambda}} {(qb)_{\lambda} (qbr^2s^2/Av^2)_{\lambda}}
    \dfrac{ (qbr/a's)_{\lambda}}{ (a's^2v^{-1})_{\lambda}}  \dfrac{ 
    (a'sv^{-1})_{\lambda}} {(qbr/a')_{\lambda}}  \\ \cdot        
  \prod_{i=1}^{n}\left\{ \dfrac{(1-brsv^{-1} t^{2-2i} q^{2\lambda_i})}
    {(1-brsv^{-1} t^{2-2i})}  \left(qt^{2i-2} \right)^{\lambda_i} \right\} 
  \prod_{1\leq i< j \leq n} \left\{ \dfrac{ (qt^{j-i})_{\lambda_i - \lambda_j}
    } { (qt^{j-i-1})_{\lambda_i - \lambda_j} } \dfrac{
    (brsv^{-1}t^{3-i-j})_{\lambda_i + \lambda_j} } { (brsv^{-1}t^{2-i-j})_{\lambda_i + 
    \lambda_j} } \right\} \\ \cdot W_{\lambda} (t, brst^{2-2n}, brsv^{-1}t^{1-n};
    q^{\nu}t^{\delta(n)};q) \, W_{\mu} (t, brsv^{-1}t^{2-2n}, bt^{1-n}; 
    q^{\lambda}t^{\delta(n)};q)   \\
= \sum_{\substack{\lambda\\ \mu \subseteq
    \lambda \subseteq \nu }} \dfrac{(r)_{\lambda}
    (Av^2/rs^2 )_{\lambda}} {(qb )_{\lambda} 
    (qbr^2s^2/Av^2)_{\lambda}} \dfrac{(qbrs/Av)_{\lambda}}{(Av/s)_{\lambda}}
    \dfrac{ (qbr/a's)_{\lambda}}{ (a's)_{\lambda}} \dfrac{ 
    (brt^{1-n})_{\lambda}}{ (qt^{n-1})_{\lambda}} \dfrac{
    (a'sv^{-1})_{\lambda}} { (qbrv/sa')_{\lambda}}\\ \cdot
 \prod_{i=1}^{n}\left\{ \dfrac{(1-brt^{2-2i} q^{2\lambda_i})}
    {(1-brt^{2-2i})}  \left(qt^{2i-2} \right)^{\lambda_i} \right\} 
  \prod_{1\leq i< j \leq n} \left\{ \dfrac{ (qt^{j-i})_{\lambda_i - \lambda_j}
    } { (qt^{j-i-1})_{\lambda_i - \lambda_j} } \dfrac{
    (brt^{3-i-j})_{\lambda_i + \lambda_j} } { (brt^{2-i-j})_{\lambda_i + 
    \lambda_j} } \right\} \\ \cdot W_{\lambda} (t, brst^{2-2n}, brt^{1-n};
    q^{\nu}t^{\delta(n)};q) \, W_{\mu} (t, brt^{2-2n},
    bt^{1-n}; q^{\lambda}t^{\delta(n)};q) 
\end{multline}
Sending $b\rightarrow 0$ and setting $A=Asr/v$ gives 
\begin{multline}
\dfrac{(v)_{\nu} } { (A)_{\mu} } 
(-1)^{|\mu|}\, A^{|\mu|}  t^{-n(\mu)} q^{n(\mu')} \\ 
\cdot \sum_{\substack{\lambda\\ \mu \subseteq
    \lambda \subseteq \nu }} \dfrac{
    (A)_{\lambda}} {  (qt^{n-1})_{\lambda} } \,      
 q^{|\lambda|} t^{2n(\lambda)}
  \prod_{1\leq i< j \leq n} \left\{ \dfrac{ (qt^{j-i})_{\lambda_i - \lambda_j}
    } { (qt^{j-i-1})_{\lambda_i - \lambda_j} } \right\} \\ 
    \cdot W^{ab}_{\lambda} (v^{-1}t^{n-1}; q^{\nu}t^{\delta(n)};t, q)\,
     W^{s}_{\mu} (q^{\lambda}t^{\delta(n)}; t, q)   \\
= \sum_{\substack{\lambda\\ \mu \subseteq
    \lambda \subseteq \nu }} \dfrac{ (-1)^{|\lambda|}\, A^{|\lambda|} v^{|\lambda|} 
    t^{-n(\lambda)} q^{n(\lambda')} } {(qt^{n-1})_{\lambda} } 
    \, q^{|\lambda|} t^{2n(\lambda)} t^{(1-n)|\lambda|} \\
\cdot  \prod_{1\leq i< j \leq n} \left\{ \dfrac{ (qt^{j-i})_{\lambda_i - \lambda_j}
    } { (qt^{j-i-1})_{\lambda_i - \lambda_j} } \right\} 
    W^{s}_{\lambda} (q^{\nu}t^{\delta(n)}; t, q) \,
    W^{s}_{\mu} (q^{\lambda}t^{\delta(n)}; t,q) 
\end{multline}
We now apply the $v=s$ and $A= -s^{-1}$ case of the last transformation identity on the sum side of~(\ref{strongci_limit}) with $z=-s^{-1}$. The result follows after some manipulations of terms and simplifications. 
\end{proof}
A multiple analogue of the identity 
\begin{equation}
2^{n-1} n = \sum_{k=1}^n 
   k \binom{n}{k} 
\end{equation}
follows from previous theorem as a special case.
\begin{cor}
For $n$-part partitions $\nu$ and $\mu$, we have
\begin{multline}
t^{-n(\mu)} q^{n(\mu')} 
\dfrac{ ( -1 )_{\nu}} { (-1 )_{\mu} } \, \left(\sum_{i=1}^{n} \dfrac{(1-q^{\nu_i})\, t^{1-i}}{1-q} \right) \\ = \sum_{e_1 \subseteq  \lambda \subseteq \nu } 
    t^{-n(\lambda)} q^{n(\lambda')}  
    \left(\sum_{i=1}^{n} \dfrac{(1-q^{\lambda_i})\, t^{1-i}}{1-q} \right)
   \binom{\nu}{\lambda}_{\!\!\!q,t}
\end{multline}
\end{cor}
 
\begin{proof}
This follows immediately by setting $\mu=e_1$ in~(\ref{doublebinom}) and using the special evaluations~(\ref{evals}) studied below.
\end{proof}

One of the most important families of discrete probability distributions is the binomial distribution whose density is given by 
\begin{equation}
f(i;n,p):= \binom{n}{i} p^i (1-p)^{n-i}
\end{equation}
where $0<p<1$ and $i=0,1,\ldots, n$. Here the density $f(i;n,p)$ denotes the probability that the ``event $i$'' occurs for fixed parameters $n$ and $p$. Changing the roles of $i$ and $n-i$ would give an equivalent definition. 

We now define an analogous probability measure on the set of all $n$-part partitions contained in $\lambda$ under the partial inclusion ordering.
\begin{definition}
Let $\lambda$ be an $n$-part partition. For any partition $\mu\subseteq\lambda$, the $qt$-binomial density function is defined by 
\begin{equation}
g(\mu;\lambda,z):= \binom{\lambda}{\mu}_{\!\!\!q,t} \, z^{|\lambda|-|\mu|} (z)_\mu  
\end{equation}
where $p\in\mathbb{R}$ with $0<p<1$. An alternative definition would be
\begin{equation}
f(\mu;\lambda,z):= t^{-2n(\mu)} q^{2n(\mu')} 
 \binom{\lambda}{\mu}_{\!\!\!q,t} \, \dfrac{(z)_\lambda} {(z)_\mu}  \,
   z^{|\mu|}
\end{equation}
\end{definition}
It needs to be verified that the density function is non--negative for any $\mu\subseteq\lambda$, and that % $g(\mu;\lambda,p)$ 
the total probability adds up to $1$ when summed over all $\mu\subseteq\lambda$. This is what we verify next.

\begin{theorem}
\label{densities}
The multiple discrete density functions $g(\mu;\lambda,z)$ and $f(\mu;\lambda,z)$ are valid densities.
\end{theorem}

\begin{proof}
Set $a=qb/r$ in~(\ref{simplifiedJs2}) and send $b\rightarrow 0$ and $s\rightarrow 0$ to get
\begin{equation}
\label{step2}
  \dfrac{(r)_\lambda }{(r/x )_\lambda}
 =\sum_{\mu \subseteq \lambda} \dfrac{(1/x )_\mu} {( r/x)_\mu}  
  (-1)^{|\mu|}\, r^{|\mu|} t^{-n(\mu)} q^{n(\mu')} 
 \binom{\lambda}{\mu}_{\!\!\!q,t}
\end{equation}
Setting $z=1/x$ now and sending $r\rightarrow\infty$ gives
\begin{equation}
  1 =\sum_{\mu \subseteq \lambda} z^{|\lambda|-|\mu|} (z)_\mu  
 \binom{\lambda}{\mu}_{\!\!\!q,t}
\end{equation}
which shows that $g(\mu;\lambda,z)$ add up to 1.
Similarly, shifting $x$ to $xr$ in~(\ref{step2}) and sending $r\rightarrow 0$ instead gives
\begin{equation}
\label{totaldensity}
1 =\sum_{\mu \subseteq \lambda} \dfrac{(z)_\lambda} {( z)_\mu}  \,
   z^{|\mu|} t^{-2n(\mu)} q^{2n(\mu')} 
 \binom{\lambda}{\mu}_{\!\!\!q,t}
\end{equation}
verifying total probability for the alternative definition $f(\mu;\lambda,z)$. 

It is obvious from the definitions~(\ref{ellipticQtPocSymbol}) and~(\ref{qtbinomcoeff}) that both densities are non--negative when $0<q,t,z<1$ and $z<t$. 
\end{proof}

For a fixed $n$-part partition $\nu$, let $S_\nu$ be the set of all partitions that are under $\nu$ with respect to the partial inclusion ordering. That is,
\[ S_\nu := \{\lambda: \lambda \subseteq \nu \} \]
It then follows from Theorem~\ref{densities} that the distribution
\begin{equation}
  F(\lambda; \nu, z) :=\sum_{\mu \subseteq \lambda} z^{|\nu|-|\mu|} (z)_\mu  
 \binom{\nu}{\mu}_{\!\!\!q,t}
\end{equation}
defines a probability measure on $S_\nu$. An alternative distribution may be defined using the $f$ density function defined above as well. The relation between this measure and the one given in~\cite{Fulman1} on the set of all partitions is to be investigated in another publication. 
   
We will, however, define a multiple analogue of another important family of density function for the Poisson distribution in this section. It will be defined as a limiting case of the binomial as in the classical case. First we give two multiple analogues of the exponential function $e^z$. 

\begin{theorem}
The following functions are multiple analogues of the exponential function $e^z$.  
\begin{multline}
 E_q(z):=(-z)_{\infty^n}
 = \sum_{\mu\in P_n} \dfrac{z^{|\mu|} q^{n(\mu')} t^{n(\mu)+(1-n)|\mu| } }  { (qt^{n-1} )_\mu} \\ \cdot \prod_{1\leq i<j \leq n} \left\{\dfrac{ (qt^{j-i})_{\mu_i-\mu_j} } {(qt^{j-i-1})_{\mu_i-\mu_j}  }
  \dfrac{(t^{j-i+1})_{\mu_i -\mu_j} } {(t^{j-i})_{\mu_i -\mu_j} } \right\}
\end{multline}
and 
\begin{multline}
e_q(z):= \dfrac{1}{(z)_{\infty^n}} 
 = \sum_{\mu \in P_n }  
 % \prod_{i=1}^n \left\{ (qt^{2i-2})^{\mu_i}\right\} 
 \dfrac{z^{|\mu|} t^{2n(\mu)+(1-n)|\mu|} } { (qt^{n-1} )_\mu} \\
\cdot \prod_{1\leq i<j \leq 
 n} \left\{\dfrac{ (qt^{j-i})_{\mu_i-\mu_j} }
 {(qt^{j-i-1})_{\mu_i-\mu_j}  } 
 \dfrac{(t^{j-i+1})_{\mu_i
-\mu_j} } {(t^{j-i})_{\mu_i -\mu_j} }
 \right\}
\end{multline}
where $\max \{\abs{zt^{(2i-n-1)} } : \, i\in[n] \} < 1$. 
\end{theorem}

\begin{proof}
Setting $\lambda=k^n$ in $qt$-binomial theorem~(\ref{qt_binom_thmAlt}), sending $k\rightarrow \infty$ by using the identity~(\ref{Wablimit}) and applying the multiple analogue of the dominated convergence theorem employed in the proof of the previous theorem, and replacing $x$ by $-z$ gives the expression for $E_q(z)$. 

Similarly, sending $x\rightarrow \infty$ in~(\ref{2phi1}), setting $\lambda=k^n$ as above and sending $k\rightarrow \infty$ using the identity~(\ref{Wslimit}), and replacing $s$ by $z$ gives $e_q(z)$. Note in the latter case that the convergence theorem requires the condition $\max \{\abs{zt^{(2i-n-1)} } \!: \! i\in[n] \} < 1$ to be satisfied. 
\end{proof}

We point out only two obvious properties of these functions. Namely, that
\[ e_q(z) \, E_q(-z) = 1 \]
and that
\[ \lim_{q\rightarrow 1}\lim_{t\rightarrow 1} E_q(z(1-q)) = e^z, \quad \mathrm{and}\quad \lim_{q\rightarrow 1}\lim_{t\rightarrow 1} e_q(z(1-q)) =e^z \]  
Both of these properties follow immediately from their definitions given in the previous theorem. 
Other properties $E_q(z)$ and $e_q(z) $ satisfy will be investigated in a future publication. 

We now give a multiple analogue of the $qt$-Poisson distribution. 
\begin{theorem}
The function
\begin{multline}
f(\mu;z):=  E_q(-z) \cdot
\dfrac{z^{|\mu|}  \,q^{2n(\mu')} t^{(1-n)|\mu| } }  { (z)_\mu(qt^{n-1} )_\mu} \\ \cdot
\prod_{1\leq i<j \leq n} \left\{\dfrac{ (qt^{j-i})_{\mu_i-\mu_j} } {(qt^{j-i-1})_{\mu_i-\mu_j}  }  \dfrac{(t^{j-i+1})_{\mu_i
-\mu_j} } {(t^{j-i})_{\mu_i -\mu_j} } \right\} 
\end{multline}
defines a valid density on the set of all partitions of length at most $n$.
\end{theorem}

\begin{proof}
Similar to the classical case, we set $\lambda=k^n$ and send $k\rightarrow\infty$ in the $qt$-binomial density $f(\mu;\lambda,z)$ using the identity~(\ref{Wslimit}) to get the $qt$-Poisson density $f(\mu;z)$ function as desired. 

That the total probability $\sum_{\mu\in P_n} f(\mu;z)$ add up to $1$, and that $f(\mu;z)$ is always non--negative follow from the construction. 
\end{proof}

Finally, we study certain special evaluations of the $qt$-binomials. It turns out that when either $\lambda$ or $\mu$ is a rectangular partition $k^n=(k,k,\ldots, k)$, then $\binom{\lambda}{\mu}_{\!q,t}$ has a closed form product representation. 

\begin{theorem}
For $n$-part partitions $\lambda$ and $\mu$, the special cases
\begin{equation}
\label{evals}
 \binom{\lambda}{k^n}_{\!\!\!q,t}, \quad  \binom{k^n}{\mu}_{\!\!\!q,t} \quad \mathrm{and}\quad \binom{\lambda}{e_1}_{\!\!\!q,t} 
\end{equation}
have closed form product representations. 
%Here
%e_1=(1,0,\ldots,0) $ denotes the $n$-part partition whose only non--zero part is 1 and is located in its %first coordinate.
\end{theorem}
 
\begin{proof}
The $\lambda=k^n$ case of identity~(\ref{commonfactor}) becomes
\begin{equation}
W_{k^n}(z;q,t,a,b) 
=\prod_{1\leq i<j\leq n}\dfrac{(qbt^{j-2i})_{2k}}{(qbt^{j-1-2i})_{2k}}
\prod_{i=1}^n\dfrac{(z_i^{-1})_k (az_i)_k}{(qbz_i)_k (qb/(az_i))_k}
\end{equation}
Set $b=as$ % to get
and send $a\rightarrow 0$ to get
\begin{equation}
W^{ab}_{k^n}(z;q,t,s) 
= \prod_{i=1}^n\dfrac{(z_i^{-1})_k }{ (qs/z_i)_k}
\end{equation}
Furthermore, multiplying both sides by $s^{nk}$ and sending $s\rightarrow \infty$ gives  
\begin{equation}
\label{Wsrect}
W^{s\uparrow}_{k^n}(z;q,t) 
= q^{-nk} \prod_{i=1}^n (q^{1-k} z_i)_k 
\end{equation}
after flipping certain factors using
\begin{equation}
\label{flips}
(a)_n = (q^{1-n}/a)_n q^{\binom{n}{2}} (-a)^n
\end{equation}
Therefore we get 
\begin{equation}
\label{rectmu_binom}
\binom{\lambda}{k^n}_{\!\!\!q,t} = \prod_{i=1}^n  \dfrac{(q^{1-k+\lambda_i}t^{n-i} )_k}{(qt^{n-i} )_{k}} 
\end{equation}

It is clear from the definition that  
\begin{equation}
W^{s\uparrow}_{\lambda/\lambda}(x; q, t) = 1
\end{equation}
and
\begin{equation}
W^{s\uparrow}_{e_1/0^n}(x; q, t) = q^{-1}(1-x)
\end{equation}
for $x\in\mathbb{C}$. It then follows from the recurrence relation~(\ref{Wsuprec}) that
\begin{equation}
W^{s\uparrow}_{e_1}(z; q, t) = \sum_{i=1}^n - q^{-1}x_i(1-x_i^{-1}t^{n-i}) = q^{-1} \left(\sum_{i=1}^n x_i + t^{i-1} \right)
\end{equation}
for $z=(x_1,\ldots,x_n)\in\mathbb{C}^n$. Multiplying by the front factors gives
\begin{equation}
\label{e1mu_binom}
\binom{\lambda}{e_1 }_{\!\!\!q,t} = \sum_{i=1}^{n} \dfrac{(1-q^{\lambda_i})\, t^{1-i}}{1-q} = \sum_{i=1}^{n} [\lambda_i]_{q} \, t^{1-i}
\end{equation}
where $[n]_{\!q} := \dfrac{1-q^n}{1-q}$ is the so--called $q$-number or $q$-bracket. 

The proof of the last result follows immediately from the identity~(\ref{Wslimit}). Setting $\lambda=k^n$ and manipulating factors gives 
\begin{multline}
\binom{k^n}{\mu}_{\!\!\!q,t} = t^{2n(\mu)+(1-n)|\mu| } \prod_{i=1}^n \dfrac{ (q^{1+k-\mu_i}t^{i-1})_{\mu_i} }  { (qt^{n-i} )_{\mu_i} } \\
\cdot \prod_{1\leq i<j \leq n} \left\{\dfrac{ (qt^{j-i})_{\mu_i-\mu_j} } {(qt^{j-i-1})_{\mu_i-\mu_j}  } \dfrac{(t^{j-i+1})_{\mu_i
-\mu_j} } {(t^{j-i})_{\mu_i -\mu_j} } \right\} 
\end{multline}
as claimed.
\end{proof} 

In the light of last theorem, we give a definition for $qt$-number $[z]_{qt}$ extending that of a $q$-number defined above. The definition will be used in the next section in the discussion of the $qt$-Stirling numbers. First, we write an extension of our $qt$-binomial coefficients.  
\begin{definition}
\label{qtbinomcoeffExt}
Let $z=(x_1,\ldots, x_n)\in\mathbb{C}^n$ and $\mu$ be $n$-part partitions. Then the extended $qt$-binomial coefficient is defined by 
\begin{equation}
\binom{z}{\mu}_{\!\!\!q,t} := \dfrac{ q^{|\mu|} t^{2n(\mu)+(1-n)|\mu| } }  { (qt^{n-1} )_\mu} \prod_{1\leq i<j \leq n} \left\{\dfrac{ (qt^{j-i})_{\mu_i-\mu_j} } {(qt^{j-i-1})_{\mu_i-\mu_j}  } \right\}  W^{s\uparrow}_\mu(q^z t^{\delta(n)}; q, t)
\end{equation}
where $q,t\in\mathbb{C}$. It should be noted that this definition makes sense even for $\mu\in\mathbb{C}^n$ by the virtue of Theorem~\ref{thm:extension} when $z$ is an $n$-part partition. 
\end{definition}

With this extension, a definition for $qt$-number may be written as follows.  
\begin{definition}
Let $z=(x_1,\ldots, x_n)\in\mathbb{C}^n$. Then 
\begin{equation}
\label{qtnumber}
[z]_{qt} =  [z; q, t] :=
\binom{z}{1^n}_{\!\!\!q,t} = \prod_{i=1}^n \dfrac{(1-q^{x_i} t^{n-i} )}{(1-qt^{n-i} )} 
\end{equation}
We also define the partition shifted generalization $[z;\mu]_{qt} = [z;\mu, q, t]$ of the $qt$-bracket in the form of a normalized $W^{s\uparrow}$ function as follows.
\begin{multline}
[z; \mu]_{qt} := 
%(-1)^{|\mu|} t^{-n(\mu)} 
q^{|\mu|+n(\mu')} 
\prod_{i=1}^n  \left\{ \dfrac{1}{(1-qt^{n-i} )^{\mu_i}} \right\} \\ \cdot
\prod_{1\leq i<j \leq n} \left\{ \dfrac{ (t^{j-i})_{\mu_i-\mu_j} } {(t^{j-i+1})_{\mu_i-\mu_j}  } \right\} \,  W^{s\uparrow}_\mu(q^z t^{\delta(n)}; q, t) 
\end{multline}
\end{definition}
The evaluation~(\ref{qtnumber}) in the definition follows from~(\ref{Wsrect}). Note also that when $z=q^x t^{\delta(n)}$ for a single variable $x\in\mathbb{C}$, the $\mu$-shifted $qt$-number $[z; \mu]_{qt}$ can be written as 
\begin{equation}
\label{qtbracket1D}
[x; \mu]_{qt} = q^{n(\mu')}
%(-1)^{|\mu|}  t^{-n(\mu)} 
(q^x;1/q,1/t)_\mu \prod_{i=1}^n \left\{ \dfrac{1}{(1-qt^{n-i} )^{\mu_i}} \right\}
\end{equation}
by using the following application of~(\ref{flips}) 
\begin{equation}
x^{|\mu|}   (x^{-1}, q, t)_\mu 
= (-1)^{|\mu|} q^{n(\mu)} t^{-n(\mu)} (x;q^{-1}, t^{-1})_{\mu}
\end{equation}
to simplify the $W^{s\uparrow}$ function.
Similarly, using~(\ref{Wsrect}) we get
\begin{equation}
[z; k^n]_{qt} := 
%(-1)^{kn}  t^{-n(k^n)}
q^{n\binom{k}{2}} 
\prod_{i=1}^n  \left\{ \dfrac{(q^{1-k} q^{x_i}t^{n-i})_k } 
{(1-qt^{n-i} )^{k}} \right\} 
\end{equation}
In particular, for $k=1$, we recover the definition~(\ref{qtnumber}) above. 

Note also that $(x;1/q,1/t)_\mu $, with the reciprocals of $q$ and $t$, corresponds to a multiple basic $qt$-analogue of the falling factorial $x_{\underline{n}}:= x (x-1) \cdots (x-(n-1))$.  

\section{Multiple basic and ordinary special numbers}
\label{mulipleSpecial}
In this section we give multiple basic and multiple ordinary analogues (or the $qt$- and $\alpha$-analogues, respectively) of several celebrated families of special numbers including Stirling numbers, Bernoulli numbers, Bell numbers, Fibonacci numbers and pentagonal numbers. The definition and some immediate properties of these numbers are studied in each case. Their fascinating deep properties such as the recurrence relations they satisfy, closed form evaluations in certain special cases and their combinatorial interpretations are carried out in other publications. 

The multiple special number sequences we introduce in this section are indexed by partitions. For each class of multiple basic $qt$-special number, the corresponding multiple ordinary $\alpha$-special number analogues can be found by setting $t=q^\alpha$ and sending $q\rightarrow 1$. 

\subsection{$qt$-Stirling numbers}
The Stirling numbers of the first kind are defined to be the coefficients of the falling factorial $x_{\underline{n}}$ in the expansion 
\begin{equation}
x_{\underline{n}} = n! \binom{x}{n}  = \sum_{k=0}^n s_1(n,k) \,x^k  
\end{equation}
The $q$-analogue of these numbers are defined in~\cite{Carlitz1} and their properties are studies in~\cite{Milne2, Gould1, Kim1, Zeng1} and others. We give the definition of the multiple $qt$-Stirling numbers generalizing the one dimensional $q$-analogues as follows.

\begin{definition}
For an $n$-part partition $\lambda$, the $qt$-Stirling numbers of the first kind $s_1(\lambda,\mu)$ are defined by 
\begin{equation}
 [x; \lambda]_{qt} 
= \sum_{\mu \subseteq \lambda} s_1(\lambda,\mu) \,  \prod_{i=1}^n \,[\,xt^{n-i}]_{qt}^{\mu_i} 
\end{equation}
where $x\in\mathbb{C}$. Consider the infinite dimensional lower triangular matrix $m$ whose entries are $s_1(\lambda,\mu)$. That is, $m_{\lambda\mu}=s_1(\lambda,\mu)$. The $qt$-Stirling numbers of the second kind $s_2(\lambda,\mu)$ are defined to be the entries of the inverse $m^{-1}$ of the matrix $m$. That is, $s_2(\lambda,\mu):=m^{-1}_{\lambda\mu}$.  
\end{definition}

We now write an explicit formula for the $qt$-Stirling numbers generalizing one dimensional $q$-analogues as follows. The proof of the theorem makes clear that the definition makes sense by presenting an inverse for the matrix $m$. 
\begin{theorem}
For $n$-part partitions $\nu$ and $\mu$, an explicit formula for the $qt$-Stirling numbers of first and second kind $s_1(\lambda,\mu)$ and $s_2(\lambda,\mu)$ are given by
\begin{equation}
s_1(\nu,\mu):=  
\dfrac{ q^{n(\nu')} t^{-2n(\mu)+(n-1)|\mu|} } 
{\prod_{i=1}^n (1-qt^{n-i} )^{\nu_i-\mu_i}} 
\cdot \!\!\!
\sum_{\mu \subseteq \lambda \subseteq \nu} \!\!\!
  u(\nu, \lambda)\, t^{(1-n)|\lambda|} \Big( \lim_{q\rightarrow 1} v(\lambda,\mu, 1/q, 1/t) \Big)
\end{equation}
and
\begin{multline}
s_2(\nu,\mu) = s_2(\nu,\mu,q,t) \\ :=
\dfrac{ q^{-n(\mu')} t^{2n(\nu)+(1-n)|\nu|} } 
{\prod_{i=1}^n (1-qt^{n-i} )^{\nu_i-\mu_i}} \cdot \!\! 
  \sum_{\mu \subseteq \lambda \subseteq \nu} \!\!
\Big( \lim_{q\rightarrow 1} u(\nu, \lambda,1/q,1/t)  \Big) t^{(n-1)|\lambda|} \, v(\lambda,\mu,q,t) 
\end{multline}
where $u(\lambda,\mu)$ and $v(\lambda,\mu)$ are defined by
\begin{equation}
u(\lambda,\mu,q, t)
:= \dfrac{ q^{|\mu|} t^{2n(\mu)}  }{(qt^{n-1})_\mu
} \prod_{1\leq i < j \leq n}
\left\{ \dfrac{ (qt^{j-i})_{\mu_i
      -\mu_j}}{(qt^{j-i-1})_{\mu_i -\mu_j} } \right\} 
 W^{s\downarrow}_\mu(q^\lambda t^{\delta(n)};q,t) 
\end{equation}
and
\begin{multline}
v(\lambda,\mu,q,t)
:= \dfrac{ (-1)^{|\mu|} q^{|\mu|+n(\mu')} t^{n(\mu)+(1-n)|\mu|} }  { (qt^{n-1} )_\mu}   \\ \cdot
\prod_{1\leq i<j \leq n} \left\{\dfrac{ (qt^{j-i})_{\mu_i-\mu_j} } {(qt^{j-i-1})_{\mu_i-\mu_j}  } \right\}  W^{s\uparrow}_\mu(q^\lambda t^{\delta(n)}; q, t) 
\end{multline}

\end{theorem}

\begin{proof}
We write the $W$-Jackson sum~\cite{CoskunG1} in the form
\begin{multline}
\label{WJackson}
W_{\lambda}(x;q,p,t,at^{-2n},bt^{-n})\\
=\dfrac{(s)_{\lambda}(as^{-1}t^{-n-1})_{\lambda}}
{(qbs^{-1}t^{-1})_{\lambda}(qbt^ns/a)_{\lambda}}
\cdot \prod_{1\leq i < j\leq n}
\left\{\dfrac{(t^{j-i+1})_{\lambda_i-\lambda_j}
(qbt^{-i-j+1})_{\lambda_i+\lambda_j}}
{(t^{j-i})_{\lambda_i-\lambda_j} (qbt^{-i-j})_{\lambda_i +
    \lambda_j}} \right\}  \\
\cdot\sum_{\mu \subseteq \lambda}
  \dfrac{(bs^{-1}t^{-n})_{\mu} (qbt^n/a)_\mu}{(qt^{n-1})_\mu
(as^{-1}t^{-n-1})_\mu} \cdot \prod_{i=1}^n \left\{
  \dfrac{(1-bs^{-1}t^{1-2i}q^{2\mu_i})}{(1-bs^{-1}t^{1-2i})}(qt^{2i-2})^{\mu_i}
\right\} \\ \cdot \prod_{1\leq i < j \leq n}
\left\{ \dfrac{(t^{j-i})_{\mu_i -\mu_j} (qt^{j-i})_{\mu_i
      -\mu_j}}{(qt^{j-i-1})_{\mu_i -\mu_j}(t^{j-i+1})_{\mu_i
      -\mu_j}} \dfrac{(bs^{-1}qt^{-i-j})_{\mu_i+\mu_j}
    (bs^{-1}t^{-i-j+2})_{\mu_i+\mu_j}} {(bs^{-1}t^{-i-j+1})_{\mu_i+\mu_j}
(qbs^{-1}t^{-i-j+1})_{\mu_i+\mu_j}} \right\} \\
\cdot W_\mu(q^\lambda t^{\delta(n)};q,t,bt^{1-2n},bs^{-1}t^{-n})\cdot
W_\mu(xs;q,t,as^{-2}t^{-2n},bs^{-1}t^{-n})
\end{multline}
Set $b=ar$ in this identity and send $a\rightarrow 0$, $r\rightarrow \infty$ using the limit rule~(\ref{LimitRule}), and set $z=xt^{\delta(n)}$ and simplify using~(\ref{Wslimit}) to get
\begin{equation}
\label{changeofbasis1}
(x;q^{-1}, t^{-1})_\lambda
  = \sum_{\mu \subseteq \lambda}
  u(\lambda, \mu) \, x^{|\mu|}  
\end{equation}
where $u(\lambda, \mu)$ is defined as above.
Similarly, replace $x$ by $x/s$ in~(\ref{WJackson}) and send $s$ to 0 to get
\begin{equation}
\label{changeofbasis2}
 x^{|\lambda|} 
 =   \sum_{\mu \subseteq \lambda}
v(\lambda, \mu)\, (x;q^{-1}, t^{-1})_{\mu}
\end{equation}
where $v(\lambda, \mu)$ is defined as in the theorem.
Then it is clear, by a change of basis argument, that 
\begin{equation}
\label{inv_uandv}
\delta_{\nu\lambda} = \sum_{\mu \subseteq \lambda \subseteq \nu}
  u(\nu, \lambda) \,  v(\lambda,\mu) 
\end{equation}
Note that the \lhs~(\ref{changeofbasis2}) does not depend on $q$ or $t$. Also,
\begin{equation}
[x t^{n-1};  \mu, 1, t] 
= \prod_{i=1}^n \,[xt^{n-i}]_{qt}^{\mu_i} 
= \prod_{i=1}^n \dfrac{(1-q^{x} t^{n-i} )^{\mu_i}}{(1-qt^{n-i} )^{\mu_i}} 
\end{equation}
for $x\in\mathbb{C}$. 
By sending $q\rightarrow 1$ in identity~(\ref{changeofbasis2}) and substituting it in~(\ref{changeofbasis1}) after replacing $x$ by $xt^{n-1}$, we get
\begin{multline}
(x;q^{-1}, t^{-1})_\lambda = \sum_{\mu \subseteq \nu} (xt^{n-1}; 1, t^{-1})_{\mu} \\ \cdot
\sum_{\mu \subseteq \lambda \subseteq \nu}
  u(\nu,\lambda) \, t^{(1-n)|\lambda|} \, \Big( \lim_{q\rightarrow 1} v(\lambda,\mu, 1/q, 1/t) \Big)
\end{multline}
This gives the explicit expression for the $qt$-Stirling numbers of the first kind after multiplying both sides by appropriate factors coming from the definition~(\ref{qtbracket1D}) of $[x; \lambda, q, t]$.
Similarly, using the inversion relation~(\ref{inv_uandv}) gives the explicit formula for the $qt$-Stirling numbers of the second kind, $s_2(\nu,\mu)$.
\end{proof}
These identities reduces in one dimensional case (i.e., $n=1$) to the formulas for the $q$-Stirling numbers given in~\cite{Kim1}. 

An obvious property that follows from the inversion relation~(\ref{inv_uandv}) in the last proof is that, considered as matrix entries, the $qt$-Stirling numbers of the first and second kind form inverse matrices. That is,
\begin{equation}
\delta_{\nu\lambda} = \sum_{\mu \subseteq \lambda \subseteq \nu}
  s_1(\nu, \lambda) \,  s_2(\lambda,\mu) = \sum_{\mu \subseteq \lambda \subseteq \nu}
  s_2(\nu, \lambda) \,  s_1(\lambda,\mu) 
\end{equation}
Leaving the investigation of other properties to a future publication, we point out another immediate property in the next theorem.
\begin{theorem}
For an $n$-part partition $\lambda$
\begin{equation}
  s_1(\lambda, \lambda) = s_2(\lambda,\lambda) = 1. 
\end{equation}
\end{theorem}

\begin{proof}
Setting $b=ast^{n-1}$ and sending $a\rightarrow 0$, and then $s\rightarrow 0$ in~(\ref{Wnormal}) gives  
\begin{multline}
W^{s\downarrow}_\lambda(q^\lambda t^{\delta(n)}; q, t) \\
=  (-1)^{|\lambda|}\, t^{-n(\lambda)} q^{-|\lambda|-n(\lambda')} 
(qt^{n-1} )_\lambda 
\cdot\prod_{1\leq i < j\leq
n} \frac{(qt^{j-i-1})_{\lambda_i-\lambda_j} }
{(qt^{j-i})_{\lambda_i-\lambda_j} }
\end{multline}
Substitute this in $u(\lambda,\lambda)$ to get
\begin{equation}
u(\lambda,\lambda)
= (-1)^{|\lambda|}\,  t^{n(\lambda)}  q^{-n(\lambda')} 
\end{equation}
Similarly, substituting~(\ref{Wsnormal}) into $v(\lambda,\lambda)$ gives
\begin{equation}
v(\lambda,\lambda)
=  (-1)^{|\lambda|} q^{n(\lambda')}  
\, t^{-n(\lambda)} 
\end{equation}
In addition, the front factors 
\begin{equation}
 \dfrac{q^{n(\nu')} t^{-2n(\mu)+(n-1)|\mu|} } %(-1)^{|\nu|-|\mu| } q^{n(\nu')} t^{-n(\nu)+n(\mu)} } 
{\prod_{i=1}^n (1-qt^{n-i} )^{\nu_i-\mu_i}}
\end{equation}
and its reciprocal reduce to $ q^{\pm n(\nu')}t^{\mp 2n(\lambda)\mp(n-1)|\lambda|} $ when $\nu=\mu$. This is precisely what we want to complete the proof, because the factors $q^{\mp n(\lambda')}$ in $u(\lambda,\lambda)$ and $v(\lambda,\lambda)$ becomes 1 in the limit as $q\rightarrow 1$.
Finally, the intermediate term and the flips $q\rightarrow 1/q$ and $t\rightarrow 1/t$ cancel with the factor $t^{-2n(\lambda)+(n-1)|\lambda|}$.
\end{proof} 
It should be pointed out again that setting $t=q^\alpha$ and sending $q\rightarrow 1$ gives the multiple ordinary $\alpha$-Stirling numbers of the first and second kind. 

\subsection{$qt$-Bernoulli numbers}
Another very important family of special numbers is the classical Bernoulli numbers $B_k$ that may be defined in terms of a recurrence relation. Namely, $B_0=1$ and $B_k$, for all integers $k\geq 1$, satisfy
\[0 = \sum_{k=0}^{n-1} \binom{n}{k} B_k \]
where $n\geq 2$.
Generalizing this approach in terms of multiple basic $qt$-binomial coefficients, we now give the definition of $qt$-Bernoulli numbers. The vector $e_1$ and the partition $\lambda^1$ are defined as above.

\begin{definition}
For an $n$-part partition $\lambda$, the $qt$-Bernoulli numbers are defined by the recurrence relation
\begin{equation}
\label{Bernoulli}
0=\sum_{\mu\subseteq \lambda}   t^{-n(\mu)} q^{n(\mu')} \binom{\lambda^1}{\mu}_{\!\!\!q,t} \beta_\mu
\end{equation}
with the initial condition that $\beta_{0^n}=\beta_{(0,0,\ldots,0)}=1$. 
\end{definition}
Solving the equation~(\ref{Bernoulli}) for $\beta_\lambda$ gives
\begin{equation}
\beta_{\lambda}= - t^{n(\lambda)} q^{-n(\lambda')}  \binom{\lambda^1}{\lambda}_{\!\!\!q,t}^{\!\!\!-1} \sum_{\mu\subset\lambda} 
 t^{-n(\mu)} q^{n(\mu')} \binom{\lambda^1}{\mu}_{\!\!\!q,t} \beta_\mu
\end{equation}
which can be used to compute $\beta_\mu$ for all partitions $\mu$. 

In the one dimensional (i.e., $n=1$) case, the definition of $qt$-Bernoulli numbers reduces essentially to $q$-Bernoulli numbers defined in~\cite{Carlitz1} and studied by~\cite{Koblitz1, Kim1, CenkciKRS1, Carlitz4} and others. In the limit as $q\rightarrow 1$ in the one dimensional case, the $qt$-Bernoulli numbers reduce to the classical Bernoulli numbers. 

We conclude this section by the following conjecture.

\begin{theorem}
In the limit as $t=q^\alpha$ and $q\rightarrow 1$, the $qt$-Bernoulli numbers $\beta_\lambda$ vanish for all $n$-part partitions $\lambda$ whose weight is $|\lambda|=2k+1$ for any integer $k\geq 1$. 
\end{theorem}

A proof of this conjecture among many other deep properties of these numbers will be investigated in a future publication. 
 
\subsection{$qt$-Bell numbers}
Another important class of special numbers intimately related to the binomial coefficients is the Bell numbers defined by 
\[B_n := \sum_{k=0}^n s_2(n,k) \]
This definition also generalize readily to the multiple case as follows. 
\begin{definition}
Let $\lambda$ be an $n$-part partition. The $qt$-Bell numbers are defined by the relation
\begin{equation}
\label{Bell}
B_{\lambda} :=\sum_{\mu\subseteq \lambda} t^{-n(\mu)} q^{n(\mu')} s_2(\lambda,\mu)
\end{equation}
where $s_2(\lambda,\mu)$ is the $qt$-Stirling numbers of the second type.
\end{definition}

The classical Bell numbers satisfy the recurrence
\[B_{n+1} = \sum_{k=0}^n \binom{n}{k} B_k \]
Among other properties, a multiple analogue of this result for $B_{\lambda}$ in terms of $\binom{\lambda}{\mu}_{\!\!\!q,t}$ will be given in another publication.

For $n=1$ this definition reduces to $q$-Bell numbers are defined and studied in~\cite{Milne1},~\cite{Wagner1} and others.  

\subsection{$qt$-Catalan numbers}
A $qt$-generalization of the Catalan numbers already defined in~\cite{GarsiaH1} and developed in~\cite{GarsiaHa1, Haiman1} and others.
% combinatorial properties are also investigated. 
These numbers are indexed by positive integers (the weights of partitions) and reduce to the $q$-Catalan numbers under the specialization of the $t$ parameter. The multiple  $qt$-Catalan numbers defined in this paper appear to be related to the definition given~\cite{GarsiaH1}, however the precise relationship is still to be given. The $qt$-Catalan numbers of this paper reduces, when $n=1$, to the one-dimensional $q$-Catalan numbers that are defined in~\cite{Carlitz3} and studied in~\cite{Andrews3, Andrews4, Stembridge1, FurlingerH1}, and others. 

The classical Catalan numbers are defined by the relation
\begin{equation}
C_n:=\dfrac{1}{n+1} \binom{2n}{n}
\label{eq:Catalan}
\end{equation}
We now give the following multiple $qt$-analogue of these numbers. 
\begin{definition}
Let $\lambda$ be an $n$-part partition. Then the $qt$-Catalan number $C_\lambda$ is defined by 
\begin{equation}
C_\lambda:=\dfrac{1}{[\lambda^1]_{\!qt}} \binom{2 \lambda}{\lambda}_{\!\!\!q,t} 
\end{equation}
where $\lambda^1 = \lambda+e_1 = (\lambda_1+1, \lambda_2, \ldots, \lambda_n)$ as before.
\end{definition}
In the special case when $\lambda=k^n$ is a rectangular partition, we get
\begin{equation}
C_{k^n}:=\dfrac{(1-qt^{n-1} )} {(1-q^{k+1} t^{n-1} )} 
\prod_{i=2}^n \dfrac{(1-qt^{n-i} )} {(1-q^{k} t^{n-i} )}
\prod_{i=1}^n \dfrac{(q^{1+k} t^{n-i})_{k} }  { (qt^{n-i} )_{k} }  \\ 
\end{equation}
using the evaluation~(\ref{rectmu_binom}) from above.
In the particular case when $k=1$, this reduces  to 
\begin{equation}
C_{1^n}:=
\prod_{i=2}^n \dfrac{(1-q^{2} t^{n-i}) }  {(1-q t^{n-i} )  } 
\end{equation}
which, as an empty product, becomes 1 when $n=1$. 

Similarly, using~(\ref{e1mu_binom}) we see that
\begin{equation}
C_{e_1}:=\dfrac{1}{[\,e_1^1]_{\!qt}} \binom{2\,e_1 }{e_1 }_{\!\!\!q,t} = 1 
\end{equation}
Other properties of $qt$-Catalan numbers will be investigated in a future publication.

\subsection{$qt$-Fibonacci numbers}
Recall that the classical Fibonacci numbers $F_n$ have the representation
\begin{equation}
F_{n+1}=\sum_{n=k+m} \binom{k}{m} 
\end{equation}
In a similar way, we extend the definition to the multiple basic case as follows.

\begin{definition}
For an $n$-part partition $\lambda$, the $qt$-Fibonacci numbers are defined by the relation 
\begin{equation}
F_{\lambda^1} :=\sum_{\nu+\mu=\lambda} q^{2n(\mu')} t^{2(n-1)|\mu|-2n(\mu)}
\binom{\nu}{\mu}_{\!\!\!qt} 
\end{equation}
where the sum runs over all $n$-part partitions $\nu\subseteq\lambda$ and $\mu\subseteq\lambda$ such that their sum equals $\lambda$. Here, the sum $\nu+\mu$ of two partitions $\nu$ and $\mu$ is defined to be the coordinate--wise summation of the two. That is,
\[\nu+\mu:=(\nu_1+\mu_1, \nu_2+\mu_2,\ldots, \nu_n+\mu_n)\]
which itself is a partition.
\end{definition}

In one dimensional case $n=1$, the $qt$-Fibonacci numbers are equivalent to the $q$-Fibonacci numbers defined in~\cite{Schur1} and studied in~\cite{Carlitz2},~\cite{Andrews2},~\cite{GarretI1}, ~\cite{Cigler1},~\cite{Cigler2},~\cite{GoytS1} and others. In particular, in the limiting case  $q\rightarrow 1$ for $n=1$, both definitions generate the classical Fibonacci sequence. 

Among a myriad of fascinating properties of the Fibonacci numbers, the fundamental recurrence relation
\[F_n = F_{n-1} + F_{n-2},\]
which is used to define the Fibonacci sequence with the initial conditions $F_0=0$ and $F_1=1$, can be generalized to $qt$-Fibonacci numbers $F_{\lambda}$ using Theorem~\ref{binomialRec}. This fundamental result and other properties will be given in a separate publication. 

We close this section by pointing out that the relation between the so--called finite form of the classical \RRis\ and the $q$-Fibonacci sequence holds true for our multiple $qt$-Fibonacci sequence. More precisely, it is known~\cite{Andrews2, Cigler1} that the $q$-Fibonacci sequence may be generated by
%The motivation behind our definition comes from the relationship between the one dimensional %$q$-Fibonacci numbers and the \RRis. More precisely, the finite form of the \RRis\ gives a formula %for the $q$-Fibonacci numbers in the form
\begin{equation}
F_{n+1}(x)=\sum_{k=0}^{ \lfloor \! \frac{n-1}{2} \!\rfloor} x^k q^{k^2-k} \binom{n-k}{k}_{\!\!q}  
\end{equation}
in the limit as $q\rightarrow 1$ for $x=1$. The finite form of \RRis\ is essentially a limiting case of the \Wt. The multiple analogue of the \Wt\ given in~\cite{Coskun1} can be written in the form   
\begin{multline}
\label{WatsonTransform}
 \dfrac{ (qb, qb/\rho_2\sigma_2)_\lambda}{(qb/\sigma_2,
  qb/\rho_2)_{\lambda} } \sum_{\mu\subseteq \lambda} q^{|\mu|}
  t^{2n(\mu)} \dfrac{ (\sigma_2, \rho_2, qb/\rho_1\sigma_1)_\mu} {(qt^{n-1},
  qb/\sigma_1, qb/\rho_1)_{\mu} } \\
  \cdot \prod_{1\leq i<j \leq n} \left\{\dfrac{
  (qt^{j-i})_{\mu_i-\mu_j} } {(qt^{j-i-1})_{\mu_i-\mu_j} } \right\} \,
  W^{ab}_\mu(q^\lambda t^{\delta(n)}; q,t, \rho_2\sigma_2t^{n-1}/qb )
   \\ = \sum_{\mu\subseteq \lambda} \left(
  \dfrac{q^3b^2}{\sigma_1\rho_1\sigma_2\rho_2} \right)^{|\mu|}
  (-1)^{|\mu|} q^{n(\mu')} t^{n(\mu)} \prod_{i=1}^n
  \left\{\dfrac{(1-bt^{2-2i}q^{2\mu_i})} {(1-b t^{2-2i})} \right\}
\\ \cdot \prod_{1\leq i<j
  \leq n} \left\{\dfrac{ (qt^{j-i})_{\mu_i-\mu_j}
  (bt^{3-i-j})_{\mu_i+\mu_j}} {(qt^{j-i-1})_{\mu_i-\mu_j} (b
  t^{2-i-j})_{\mu_i+\mu_j} } \right\}\, W^a_\mu(q^\lambda t^{\delta(n)};
  q,t, bt^{1-n}) \\ \cdot \dfrac{(bt^{1-n}, \sigma_2, \rho_2, \sigma_1,
  \rho_1)_{\mu}}{(qt^{n-1}, qb/\sigma_1, qb/\rho_1, qb/\sigma_2,
  qb/\rho_2)_{\mu}} 
\end{multline}
The limiting case of this identity, as $\sigma_1,\rho_1,\sigma_2,\rho_2 \rightarrow \infty$, becomes 
\begin{multline}
(qb)_\lambda \sum_{\mu\subseteq \lambda} q^{|\mu|+2n(\mu')} t^{-2n(\mu)} b^{|\mu|} \binom{\lambda}{\mu}_{\!\!\!q,t} 
   \\ = 
 \sum_{\mu\subseteq \lambda} (-1)^{|\mu|} q^{3|\mu|+5n(\mu')} t^{-3n(\mu)} b^{2|\mu|}
\dfrac{   (bt^{1-n} )_{\mu}}{(qt^{n-1} )_{\mu} } \prod_{i=1}^n
  \left\{\dfrac{(1-bt^{2-2i}q^{2\mu_i})} {(1-b t^{2-2i})} \right\}
\\ \cdot \prod_{1\leq i<j
  \leq n} \left\{\dfrac{ (qt^{j-i})_{\mu_i-\mu_j}
  (bt^{3-i-j})_{\mu_i+\mu_j}} {(qt^{j-i-1})_{\mu_i-\mu_j} (b
  t^{2-i-j})_{\mu_i+\mu_j} } \right\}\, W^a_\mu(q^\lambda t^{\delta(n)};
  q,t, bt^{1-n}) 
\end{multline}
which is a multiple analogue of the finite (terminating) version of the multiple \RRis\ studied in~\cite{Coskun1}. The summand on the \lhs\ becomes identical to that in our definition of $qt$-Fibonacci numbers on making the substititution $b=q^{-1}t^{2(n-1)}$ except that the only partitions $\nu$ and $\mu$ that enter the sum in the definition are the ones that add up to $\lambda$. 

\section{Conclusion}
The classical sequences of special numbers have fascinated mathematicians for centuries, and great many properties of these numbers have been studied. Several one or two parameter extensions of these numbers have been developed as well. However, these extensions were mostly one dimensional, and multiple analogues of these numbers have been missing. This paper takes an important step and presents multiple (basic and ordinary) analogues of several important classes of special numbers and develops some of their basic properties. Needless to say, the greatly rich properties of these numbers can not be fully covered in a single article. In addition, no attempt is made to include other families of related special numbers in this paper. These ideas will be pursued in future publications.

\end{document}